
%
%

\overfullrule=0pt

\def\bkt{\it} 

\font\headerfont=cmr10

\font\titlefont=cmbx12 scaled\magstephalf
\font\sectionfont=cmbx12
\font\foott=cmtt9
\font\footfont=cmr8
\font\footfonts=cmr6

\font\scriptit=cmti10 at 7pt
\font\scriptsl=cmsl10 at 7pt
\scriptfont\itfam=\scriptit
\scriptfont\slfam=\scriptsl

\newif\ifblackboardbold
\blackboardboldtrue

\newfam\bboldfam
\ifblackboardbold
\font\tenbbold=msbm10
\font\sevenbbold=msbm7
\font\fivebbold=msbm5
\textfont\bboldfam=\tenbbold
\scriptfont\bboldfam=\sevenbbold
\scriptscriptfont\bboldfam=\fivebbold
\def\bbold{\fam\bboldfam\tenbbold}
\else
\def\bbold{\bf}
\fi

\newcount\amsfamcount 
\newcount\classcount   
\newcount\positioncount
\newcount\codecount
\newcount\nct             
\def\newsymbol#1#2#3#4#5{               
\nct="#2                                  
\ifnum\nct=1 \amsfamcount=\msamfam\else   
\ifnum\nct=2 \amsfamcount=\msbmfam\else   
\ifnum\nct=3 \amsfamcount=\eufmfam
\fi\fi\fi
\multiply\amsfamcount by "100           
\classcount="#3                 
\multiply\classcount by "1000           
\positioncount="#4#5            
\codecount=\classcount                  
\advance\codecount by \amsfamcount      
\advance\codecount by \positioncount
\mathchardef#1=\codecount}              


\font\Arm=cmr9
\font\Ai=cmmi9
\font\Asy=cmsy9
\font\Abf=cmbx9
\font\Brm=cmr7
\font\Bi=cmmi7
\font\Bsy=cmsy7
\font\Bbf=cmbx7
\font\Crm=cmr6
\font\Ci=cmmi6
\font\Csy=cmsy6
\font\Cbf=cmbx6

\ifblackboardbold
\font\Abbold=msbm10 at 9pt
\font\Bbbold=msbm7
\font\Cbbold=msbm5 at 6pt
\fi


\relax



\def\small{%
\textfont0=\Arm \scriptfont0=\Brm \scriptscriptfont0=\Crm
\textfont1=\Ai \scriptfont1=\Bi \scriptscriptfont1=\Ci
\textfont2=\Asy \scriptfont2=\Bsy \scriptscriptfont2=\Csy
\textfont\bffam=\Abf \scriptfont\bffam=\Bbf \scriptscriptfont\bffam=\Cbf
\def\rm{\fam0\Arm}\def\mit{\fam1}\def\oldstyle{\fam1\Ai}%
\def\bf{\fam\bffam\Abf}%
\ifblackboardbold
\textfont\bboldfam=\Abbold
\scriptfont\bboldfam=\Bbbold
\scriptscriptfont\bboldfam=\Cbbold
\def\bbold{\fam\bboldfam\Abbold}%
\fi
\rm
}

\font\tenmsb=msbm10
\font\sevenmsb=msbm10 at 7pt
\font\fivemsb=msbm10 at 5pt
\newfam\msbfam
\textfont\msbfam=\tenmsb
\scriptfont\msbfam=\sevenmsb
\scriptscriptfont\msbfam=\fivemsb

\def\hexnumber#1{\ifcase#1 0\or1\or2\or3\or4\or5\or6\or7\or8\or9\or
	A\or B\or C\or D\or E\or F\fi}

\mathchardef\subsetneq="2\hexnumber\msbfam28


\ifx\begin\undefined\else\global\advance\srcdepth by
1\expandafter\endinput\fi

\def\begin{}
\newcount\srcdepth
\srcdepth=1

\outer\def\bye{\global\advance\srcdepth by -1
  \ifnum\srcdepth=0
    \def\endcmd{\vfill\supereject\nopagenumbers\par\vfill\supereject\end}
  \else\def\endcmd{}\fi
  \endcmd
}



\def\initialize#1#2#3#4#5#6{
  \ifnum\srcdepth=1
  \magnification=#1
  \hsize = #2
  \vsize = #3
  \hoffset=#4
  \advance\hoffset by -\hsize
  \divide\hoffset by 2
  \advance\hoffset by -1truein
  \voffset=#5
  \advance\voffset by -\vsize
  \divide\voffset by 2
  \advance\voffset by -1truein
  \advance\voffset by #6
  \baselineskip=13pt
  \emergencystretch = 0.05\hsize
  \fi
}

\def\irenasize{\initialize{1195}
  {6.5truein}{8.8truein}{8.5truein}{11truein}{0.2truein}}



%






\newlinechar=`@
\def\forwardmsg#1#2#3{\immediate\write16{@*!*!*!* forward reference should
be: @\noexpand\forward{#1}{#2}{#3}@}}
\def\nodefmsg#1{\immediate\write16{@*!*!*!* #1 is an undefined reference@}}

\def\forwardsub#1#2{\def\newref{{#2}{#1}}}

\def\forward#1#2#3{%
\expandafter\expandafter\expandafter\forwardsub\expandafter{#3}{#2}
\expandafter\ifx\csname#1\endcsname\relax\else%
\expandafter\ifx\csname#1\endcsname\newref\else%
\expandafter\ifx\csname#2\endcsname\relax\else%
\forwardmsg{#1}{#2}{#3}\fi\fi\fi%
\expandafter\let\csname#1\endcsname\newref}

\def\firstarg#1{\expandafter\argone #1}\def\argone#1#2{#1}
\def\secondarg#1{\expandafter\argtwo #1}\def\argtwo#1#2{#2}

\def\ref#1{\expandafter\ifx\csname#1\endcsname\relax%
  {\nodefmsg{#1}\bf`#1'}\else%
  \expandafter\firstarg\csname#1\endcsname%
  ~\htmllocref{#1}{\expandafter\secondarg\csname#1\endcsname}\fi}

\def\refscor#1{\expandafter\ifx\csname#1\endcsname\relax
  {\nodefmsg{#1}\bf`#1'}\else
  Corollaries~\htmllocref{#1}{\expandafter\secondarg\csname#1\endcsname}\fi}

\def\refs#1{\expandafter\ifx\csname#1\endcsname\relax
  {\nodefmsg{#1}\bf`#1'}\else
  \expandafter\firstarg\csname #1\endcsname
  s~\htmllocref{#1}{\expandafter\secondarg\csname#1\endcsname}\fi}

\def\refn#1{\expandafter\ifx\csname#1\endcsname\relax
  {\nodefmsg{#1}\bf`#1'}\else
  \htmllocref{#1}{\expandafter\secondarg\csname #1\endcsname}\fi}

\def\pageref#1{\expandafter\ifx\csname#1\endcsname\relax
  {\nodefmsg{#1}\bf`#1'}\else
  \expandafter\firstarg\csname#1\endcsname%
  ~\htmllocref{#1}{\expandafter\secondarg\csname#1\endcsname}\fi}

\def\pagerefs#1{\expandafter\ifx\csname#1\endcsname\relax
  {\nodefmsg{#1}\bf`#1'}\else
  \expandafter\firstarg\csname#1\endcsname%
  s~\htmllocref{#1}{\expandafter\secondarg\csname#1\endcsname}\fi}

\def\pagerefn#1{\expandafter\ifx\csname#1\endcsname\relax
  {\nodefmsg{#1}\bf`#1'}\else
  \htmllocref{#1}{\expandafter\secondarg\csname#1\endcsname}\fi}

\def\slidepageref#1{\expandafter\ifx\csname#1\endcsname\relax%
  {\nodefmsg{#1}\bf`#1'}\else%
  \htmllocref{#1}{\expandafter\csname#1\endcsname}{1}\fi}

\def\slidepageref#1{\expandafter\ifx\csname#1\endcsname\relax%
  {\nodefmsg{#1}\bf`#1'}\else%
  \htmllocref{#1}{\expandafter\csname#1\endcsname}{ }\fi}

\def\blueslidepageref#1{\expandafter\ifx\csname#1\endcsname\relax%
  {\nodefmsg{#1}\bf`#1'}\else%
  \htmllocref{#1}{\expandafter\csname#1\endcsname}{2}\fi}

\def\blueslidepageref#1{\expandafter\ifx\csname#1\endcsname\relax%
  {\nodefmsg{#1}\bf`#1'}\else%
  \htmllocref{#1}{\expandafter\csname#1\endcsname}{ }\fi}

\def\slidepagelabel#1#2{%
   \write\isauxout{\noexpand\forward{\noexpand#1}{}{#2}}%
   \bgroup\htmlanchor{#1}{\makeref{\noexpand#1}{}{{#2}}}\egroup%
 }








\newif\ifhypers  
\hypersfalse


\edef\freehash{\catcode`\noexpand\#=\the\catcode`\#}%
\catcode`\#=12
\freehash
\let\freehash=\relax
\ifhypers\fi
\def\puthtml#1{\ifhypers\fi}
\def\htmlanchor#1#2{\puthtml{<a name="#1">}#2\puthtml{</a>}}
\def\@pdfm@mark#1{}
\def\setlink#1{\colored{\linkcolor}{#1}}%
\def\setlink#1{\ifdraft\Purple{#1}\else{#1}\fi}%

%

%
\def\htmllocref#1#2{\ifhypers\leavevmode\fi\setlink{#2}\ifhypers\fi\relax}%
%
%
%
%
\def\Acrobatmenu#1#2{%
  \@pdfm@mark{%
    bann <<
      /Type /Annot
      /Subtype /Link
      /A <<
        /S /Named
        /N /#1
      >>
      /Border [\@pdfborder]
      /C [\@menubordercolor]
    >>%
   }%
  \Hy@colorlink{\@menucolor}#2\Hy@endcolorlink
  \@pdfm@mark{eann}%
}
\def\@pdfborder{0 0 1}
\def\@menubordercolor{1 0 0}
\def\@menucolor{red}

\def\ifempty#1#2\endB{\ifx#1\endA}
\def\makeref#1#2#3{\ifempty#1\endA\endB\else\forward{#1}{#2}{#3}\fi\unskip}


\def\crosshairs#1#2{
  \dimen1=.002\drawx
  \dimen2=.002\drawy
  \ifdim\dimen1<\dimen2\dimen3\dimen1\else\dimen3\dimen2\fi
  \setbox1=\vclap{\vline{2\dimen3}}
  \setbox2=\clap{\hline{2\dimen3}}
  \setbox3=\hstutter{\kern\dimen1\box1}{4}
  \setbox4=\vstutter{\kern\dimen2\box2}{4}
  \setbox1=\vclap{\vline{4\dimen3}}
  \setbox2=\clap{\hline{4\dimen3}}
  \setbox5=\clap{\copy1\hstutter{\box3\kern\dimen1\box1}{6}}
  \setbox6=\vclap{\copy2\vstutter{\box4\kern\dimen2\box2}{6}}
  \setbox1=\vbox{\offinterlineskip\box5\box6}
  \smash{\vbox to #2{\hbox to #1{\hss\box1}\vss}}}

\def\boxgrid#1{\rlap{\vbox{\offinterlineskip
  \setbox0=\hline{\wd#1}
  \setbox1=\vline{\ht#1}
  \smash{\vbox to \ht#1{\offinterlineskip\copy0\vfil\box0}}
  \smash{\vbox{\hbox to \wd#1{\copy1\hfil\box1}}}}}}

\def\drawgrid#1#2{\vbox{\offinterlineskip
  \dimen0=\drawx
  \dimen1=\drawy
  \divide\dimen0 by 10
  \divide\dimen1 by 10
  \setbox0=\hline\drawx
  \setbox1=\vline\drawy
  \smash{\vbox{\offinterlineskip
    \copy0\vstutter{\kern\dimen1\box0}{10}}}
  \smash{\hbox{\copy1\hstutter{\kern\dimen0\box1}{10}}}}}

\long\def\boxit#1#2#3{\hbox{\vrule
  \vtop{%
    \vbox{\hsize=#2\hrule\kern#1%
      \hbox{\kern#1#3\kern#1}}%
    \kern#1\hrule}%
    \vrule}}
\long\def\boxitv#1#2#3{\boxit{#1}{#2}{{\hbox to #2{\vrule%
  \vbox{\hsize=#2\hrule\kern#1%
      \hbox{\kern#1#3\kern#1}}%
    \kern#1\hrule}%
    \vrule}}}

\newdimen\boxingdimen
\long\def\boxing#1{\boxingdimen=\hsize\advance\boxingdimen by -2ex%
\vskip0.5cm\boxit{1ex}{\boxingdimen}{\vbox{#1}}\vskip0.5cm}

\newdimen\boxrulethickness \boxrulethickness=.4pt
\newdimen\boxedhsize
\newbox\textbox 
\newdimen\originalbaseline 
\newdimen\hborderspace\newdimen\vborderspace
\hborderspace=3pt \vborderspace=3pt

\def\preparerulebox#1#2{\setbox\textbox=\hbox{#2}%
   \originalbaseline=\vborderspace
   \advance\originalbaseline\boxrulethickness
   \advance\originalbaseline\dp\textbox
   \def\Borderbox{\vbox{\hrule height\boxrulethickness
     \hbox{\vrule width\boxrulethickness\hskip\hborderspace
     \vbox{\vskip\vborderspace\relax#1{#2}\vskip\vborderspace}%
     \hskip\hborderspace\vrule width\boxrulethickness}%
     \hrule height\boxrulethickness}}}

\def\hrulebox#1{\hbox{\preparerulebox{\hbox}{#1}%
   \lower\originalbaseline\Borderbox}}
\def\vrulebox#1#2{\vbox{\preparerulebox{\vbox}{\hsize#1#2}\Borderbox}}
\def\parrulebox#1\par{\boxedhsize=\hsize\advance\boxedhsize
   -2\boxrulethickness \advance\boxedhsize-2\hborderspace
   \vskip3\parskip\vrulebox{\boxedhsize}{#1}\vskip\parskip\par}
\def\sparrulebox#1\par{\vskip1truecm\boxedhsize=\hsize\advance\boxedhsize
   -2\boxrulethickness \advance\boxedhsize-2\hborderspace
   \vskip3\parskip\vrulebox{\boxedhsize}{#1}\vskip.5truecm\par}



%

\newcount\sectno \sectno=0
\newcount\appno \appno=64 
\newcount\thmno \thmno=0
\newcount\randomct \newcount\rrandomct
\newif\ifsect \secttrue
\newif\ifapp \appfalse
\newif\ifnumct \numcttrue
\newif\ifchap \chapfalse
\newif\ifdraft \draftfalse
\newif\ifcolor\colorfalse


\def\widow#1{\vskip 0pt plus#1\vsize\goodbreak\vskip 0pt plus-#1\vsize}


\def\stdskip{\vskip 9pt plus3pt minus 3pt}
\def\stdbreak{\par\removelastskip\penalty-100\stdskip}
\def\halfbreak{\vskip 0.6ex\penalty-100}

\def\proof{\stdbreak\noindent%
  \ifdraft\let\proofcolor=\Purple\else\let\proofcolor=\Black\fi%
  \proofcolor{{\sl Proof. }}}

\def\proofone{\stdbreak\noindent%
  \ifdraft\let\proofcolor=\Purple\else\let\proofcolor=\Black\fi%
  \proofcolor{{\sl Proof \#1:\ \ }}}
\def\prooftwo{\stdbreak\noindent%
  \ifdraft\let\proofcolor=\Purple\else\let\proofcolor=\Black\fi%
  \proofcolor{{\sl Proof \#2:\ \ }}}
\def\proofthree{\stdbreak\noindent%
  \ifdraft\let\proofcolor=\Purple\else\let\proofcolor=\Black\fi%
  \proofcolor{{\sl Proof \#3:\ \ }}}
\def\proofof#1{\stdbreak\noindent%
  \ifdraft\let\proofcolor=\Purple\else\let\proofcolor=\Black\fi%
  \proofcolor{{\sl Proof of #1:\ }}}
\def\claim{\stdbreak\noindent%
  \ifdraft\let\proofcolor=\Purple\else\let\proofcolor=\Black\fi%
  \proofcolor{{\sl Claim:\ }}}
\def\proofclaim{\stdbreak\noindent%
  \ifdraft\let\proofcolor=\Purple\else\let\proofcolor=\Black\fi%
  \proofcolor{{\sl Proof of the claim:\ }}}

\newif\ifnumberbibl \numberbibltrue
\newdimen\labelwidth

\def\references{\bgroup
  \edef\numtoks{}%
  \global\thmno=0
  \setbox1=\hbox{[999]} 
  \labelwidth=\wd1 \advance\labelwidth by 2.5em
  \ifnumberbibl\advance\labelwidth by -2em\fi
  \parindent=\labelwidth \advance\parindent by 0.5em 
  \removelastskip
  \widow{0.1}
  \vskip 24pt plus 6pt minus 6 pt
  \frenchspacing
  \immediate\write\isauxout{\noexpand\forward{Bibliography}{}{\the\pageno}}%
  \immediate\write\iscontout{\noexpand\contnosectlist{}{Bibliography}{\the\pageno}}%
  \ifdraft\let\refcolor=\Maroon\else\let\refcolor=\Black\fi%
  \leftline{\sectionfont\refcolor{References}}
  \ifhypers%
     \global\thmno=0\relax%
     \ifsect%
	\global\advance\sectno by 1%
	\edef\numtoks{\number\sectno}%
  	\hbox{}%
     \else%
	\edef\numtoks{}%
  	\hbox{}%
     \fi%
  \fi%
  \nobreak\stdskip\noindent}%
\def\endreferences{\nonfrenchspacing\egroup}

\def\referencesn{\bgroup
  \edef\numtoks{}%
  \global\thmno=0
  \setbox1=\hbox{[999]} 
  \labelwidth=\wd1 \advance\labelwidth by 2.5em
  \ifnumberbibl\advance\labelwidth by -2em\fi
  \parindent=\labelwidth \advance\parindent by 0.5em 
  \removelastskip
  \widow{0.1}
  \vskip 24pt plus 6pt minus 6 pt
  \frenchspacing
  \immediate\write\isauxout{\noexpand\forward{Bibliography}{}{\the\pageno}}%
  \immediate\write\iscontout{\noexpand\contnosectlist{}{Bibliography}{\the\pageno}}%
  \ifdraft\let\refcolor=\Maroon\else\let\refcolor=\Black\fi%
  \ifhypers%
     \global\thmno=0\relax%
     \ifsect%
	\global\advance\sectno by 1%
	\edef\numtoks{\number\sectno}%
  	\hbox{}%
     \else%
	\edef\numtoks{}%
  	\hbox{}%
     \fi%
  \fi%
  \nobreak\stdskip\noindent}%

\def\bitem#1{\global\advance\thmno by 1%
  \ifdraft\let\itemcolor=\Purple\else\let\itemcolor=\Black\fi%
  \outer\expandafter\def\csname#1\endcsname{\the\thmno}%
  \edef\numtoks{\number\thmno}%
  \ifhypers\htmlanchor{#1}{\makeref{\noexpand#1}{REF}{\numtoks}}\fi%
  \ifnumberbibl
    \immediate\write\isauxout{\noexpand\forward{\noexpand#1}{}{\the\thmno}}
    \item{\hbox to \labelwidth{\itemcolor{\hfil\the\thmno.\ \ }}}
  \else
    \immediate\write\isauxout{\noexpand\expandafter\noexpand\gdef\noexpand\csname #1\noexpand\endcsname{#1}}
    \item{\hbox to \labelwidth{\itemcolor{\hfil#1\ \ }}}
  \fi}

\newcount\chapno
\newcount\chappage
\chapno=0

\newtoks\leftpagehead
\newtoks\rightpagehead
\newtoks\chapname
\newtoks\chapkludge
\footline={\bgroup{\ifdraft\ifnum\pageno>5\hfill\hbox{\headerfont\today}\fi\fi}\egroup}


\def\chapter#1{\vfill\eject%
  \chaptrue%
  \global\advance\chapno by 1%
  \edef\chaptitl{#1}%
  \edef\numtoks{\number\chapno}%
  \edef\chaptitle{Chapter \numtoks:  \chaptitl}%
  \edef\secttitle{Chapter \numtoks:  \chaptitl}%
  \edef\chapname{Chapter}%
  \edef\chapkludge{1}%
  \dochapter{#1}}

\def\chapteroddpage#1{\vfill\eject \ifodd\pageno \else\blankpage\fi%
\chapter{#1}}

\def\dochapter#1{%
  \sectfalse
  \chappage=\the\pageno%
  \global\sectno=0\relax%
  \global\thmno=0\relax%
  \global\randomct=0\relax%
  \ifdraft\let\chapcolor=\Apricot\else\let\chapcolor=\Black\fi%
  \message{#1}%
  \ifhypers\hbox{}\fi%
\footline={}%
\leftpagehead={\Black{\headerfont\the\pageno\hskip3em\chaptitle\hfill}}
\rightpagehead={\Black{\headerfont\hfill\secttitle\botmark\kern3em\the\pageno}}%
\headline={{\ifnum\pageno>\chappage%
    \ifodd\pageno%
        \the\rightpagehead%
   \else%
	\the\leftpagehead\fi\fi%
   }}%
    \futurelet\testchar\maybeoptionchapter}

\def\maybeoptionchapter{\ifx[\testchar\let\next\optionchapter%
	\else\let\next\nooptionchapter\fi\next}

\def\optionchapter[#1]{%
  \immediate\write\isauxout{\noexpand\forward{\noexpand#1}{\chapname}{\numtoks}}%
  \ifcase\chapkludge%
  \leftline{{\draftlabel{#1}\titlefont \chapcolor{\chaptitl}}}%
  \immediate\write\iscontout{\noexpand\basicchaplist{\noexpand#1}{\chaptitl}{\the\pageno}}%
  \or%
  \leftline{{\draftlabel{#1}\titlefont \chapname \ \numtoks: \chapcolor{\chaptitl}}}%
  \immediate\write\iscontout{\noexpand\chaplist{\noexpand#1}{\chaptitl}{\the\pageno}}%
  \fi%
  \htmlanchor{#1}{\makeref{#1}{\chapname}{\numtoks}}%
  \nobreak\vskip 2cm}

\def\nooptionchapter{%
  \leftline{{\titlefont \chapname \ \numtoks: \chapcolor{\chaptitl}}}%
  \immediate\write\iscontout{\noexpand\chaplist{0}{\chaptitle}{\the\pageno}}%
  \nobreak\vskip 2cm}

\outer\def\nosection#1{%
  \appfalse\sectfalse\chapfalse%
  \removelastskip%
  \global\thmno=0\relax%
  \global\randomct=0\relax%
  \chappage=\the\pageno%
  \edef\numtoks{}%
  \def\secttitl{#1}%
  \vskip 24pt plus 6pt minus 6 pt%
  \ifdraft\let\sectcolor=\OliveGreen\else\let\sectcolor=\Black\fi%
  \widow{.03}%
  \message{#1}%
  \ifhypers\hbox{}\fi%
  \def\secttitle{\secttitl}%
  \immediate\write\isauxout{\noexpand\forward{\noexpand#1}{}{\numtoks}}%
  {\noindent{\sectionfont\sectcolor{\secttitl}}}%
  \immediate\write\iscontout{\noexpand\contnosectlist{\noexpand#1}{\noexpand\noindent\ \ \secttitl}{\the\pageno}}%
  \htmlanchor{#1}{\makeref{#1}{}{\numtoks}}%
  \nobreak\vskip 4ex}

\outer\def\section#1{%
  \widow{0.1}%
  \global\subsectno=0%
  \appfalse\secttrue%
  \removelastskip%
  \global\advance\sectno by 1%
  \global\thmno=0\relax%
  \global\randomct=0\relax%
  \edef\numtoks{\ifchap\number\chapno.\fi\number\sectno}%
  \edef\secttitl{#1}%
  \edef\secttitle{Section \numtoks: \secttitl}%
  \rightpagehead={\Black{\headerfont\hfill\secttitle\kern3em\the\pageno}}%
  \vskip 24pt plus 6pt minus 6 pt%
  \ifdraft\let\sectcolor=\OliveGreen\else\let\sectcolor=\Black\fi%
  \message{#1}%
  \ifhypers\hbox{}\fi%
    \futurelet\testchar\maybeoptionsection}

\def\maybeoptionsection{\ifx[\testchar\let\next\optionsection%
	\else\let\next\nooptionsection\fi\next}

\def\optionsection[#1]{%
  \immediate\write\isauxout{\noexpand\forward{\noexpand#1}{Section}{\numtoks}}%
  {\noindent{\draftlabel{#1}\sectionfont\sectcolor{\numtoks}\quad \sectcolor{\secttitl}}}%
  \immediate\write\iscontout{\noexpand\contlist{\noexpand#1}{\secttitl}{\the\pageno}}%
  \htmlanchor{#1}{\makeref{#1}{Section}{\numtoks}}%
  \nobreak\vskip 4ex}

\def\nooptionsection{%
  {\noindent{\sectionfont\sectcolor{\numtoks}\quad \sectcolor{\secttitl}}}%
  \write\iscontout{\noexpand\contlist{0}{\secttitle}{\the\pageno}}%
  \nobreak\vskip 4ex}

\newcount\subsectno
\outer\def\subsection#1{%
  \advance\subsectno by 1%
  \global\randomct=0\relax%
  \vskip 10pt plus 6pt minus 6 pt%
  \widow{.02}%
  \message{#1}%
  \noindent \sectcolor{$\underline{\hbox{\numtoks.\number\subsectno\quad #1}}$}%
  \nobreak\vskip 1ex}

\outer\def\appendix#1{%
  \global\subsectno=0%
  \vfill\eject%
  \chappage=\the\pageno%
  \apptrue\sectfalse%
  \removelastskip%
  \global\advance\appno by 1%
  \global\thmno=0\relax%
  \global\randomct=0\relax%
  \edef\numtoks{\char\number\appno}%
  \def\secttitl{#1}%
  \def\chaptitl{#1}%
  \edef\chaptitle{Appendix \numtoks:  \chaptitl}%
  \edef\secttitle{Appendix \numtoks:  \chaptitl}%
  \vskip 24pt plus 6pt minus 6 pt%
  \ifdraft\let\appendixcolor=\Blue\else\let\appendixcolor=\Black\fi%
  \widow{.03}%
  \message{#1}%
  \ifhypers\hbox{}\fi%
    \futurelet\testchar\maybeoptionappendix}

\def\maybeoptionappendix{\ifx[\testchar\let\next\optionappendix%
	\else\let\next\nooptionappendix\fi\next}

\def\optionappendix[#1]{%
  \immediate\write\isauxout{\noexpand\forward{\noexpand#1}{Appendix}{\numtoks}}%
  {\noindent{\draftlabel{#1}\titlefont Appendix \appendixcolor{\numtoks}.\quad \appendixcolor{\chaptitl}}}%
  \immediate\write\iscontout{\noexpand\chaplist{\noexpand#1}{\secttitl}{\the\pageno}}%
  \immediate\htmlanchor{#1}{\makeref{#1}{Appendix}{\numtoks}}%
  \immediate\def\chaptitle{Appendix \numtoks: \chaptitl}%
  \immediate\def\secttitle{\chaptitle}%
  \nobreak\vskip 2cm}

\def\nooptionappendix{%
  {\noindent{\titlefont Appendix \appendixcolor{\numtoks}.\quad \appendixcolor{\chaptitl}}}%
  \immediate\write\iscontout{\noexpand\chaplist{0}{\secttitle}{\the\pageno}}%
  \def\chaptitle{Appendix \numtoks: \chaptitl}%
  \def\secttitle{\chaptitle}%
  \nobreak\vskip 2cm}

\outer\def\oldappendix#1#2{%
  \removelastskip%
  \global\advance\appno by 1%
  \global\thmno=0\relax%
  \global\randomct=0\relax%
  \edef\numtoks{Appendix\ \char\number\appno}%
  \vskip 24pt plus 6pt minus 6 pt%
  \widow{.03}%
  \leftline{\draftlabel{#1}\sectionfont\appendixcolor{\numtoks}\quad \appendixcolor{#2}}%
  \immediate\write\isauxout{\noexpand\forward{\noexpand#1}{}{\numtoks}}%
  \immediate\write\iscontout{\noexpand\contnosectlist{}{\numtoks:  #2}{\the\pageno}}%
  \htmlanchor{#1}{\makeref{#1}{Section}{\numtoks}}%
  \message{#2}%
  \ifhypers\hbox{}\fi%
  \apptrue\sectfalse%
  \nobreak\vskip 3pt}


\newif\ifnewgroup\newgroupfalse

\def\proclamation#1#2#3{
  \outer\expandafter\def\csname#1\endcsname{%
  \global\newgroupfalse%
  \ifnum#3<5\global\newgrouptrue\fi%
  \ifnum#3<1\global\newgroupfalse\fi%
  \ifdraft\let\proclaimcolor=\Fuchsia\else\let\proclaimcolor=\Black\fi%
  \gdef\Prnm{#2}%
  \global\advance\thmno by 1%
  \global\randomct=0%
  \ifcase#3
	\stdbreak \or 
	\stdbreak \or 
	\stdbreak \or 
	\stdbreak \or 
	\halfbreak \or 
	\halfbreak \or 
	\halfbreak \or 
	\halfbreak \or 
	\halfbreak \or 
	\halfbreak \or 
	\halfbreak \or 
	\else \stdbreak \fi%
  \edef\proctoks{\ifchap\the\chapno.\fi\ifsect\the\sectno.\fi\the\thmno\ifnum#3=2'\fi\ifnum#3=3''\fi}%
  \widow{0.05}%
  \ifnumct\noindent{\proclaimcolor{%
	\ifnum#3=8*\fi%
	\ifnum#3=9*\fi%
	\ifnum#3=7\hbox to 1ex{$\dag$}\fi%
	\ifnum#3=10\hbox to 0.5ex{$\dag$}\fi%
	\ifnum#3=-1{\bf Important}\ \lowercase\fi%
	\ifnum#3=6\else\ifnum#3=9\else\ifnum#3=10\else{\bf #2}\fi\fi\fi%
	\ifnum#3=6\else\ifnum#3=9\else\ \fi\fi%
	\bf \proctoks}\enspace}%
  \else\noindent{\proclaimcolor{\bf #2}\enspace}%
  \fi%
  \futurelet\testchar\maybeoptionproclaim}}

\def\maybeoptionproclaim{\ifx[\testchar\let\next\optionproclaim%
	\else\let\next\nooptionproclaim\fi\next}

\def\optionproclaim[#1]{%
  \bgroup\htmlanchor{#1}{\makeref{\noexpand#1}{\Prnm}{\proctoks}}\egroup%
  \immediate\write\isauxout{\noexpand\forward{\noexpand#1}{\Prnm}{\proctoks}}%
  \draftlabel{#1}%
  \ifnewgroup\bgroup\sl\fi}
\def\nooptionproclaim{\ifnewgroup\bgroup\sl\fi}
\def\endb{\par\stdbreak\egroup\newgroupfalse\randomct=0}

\def\pagelabel#1{%
   \ \unskip%
   \write\isauxout{\noexpand\forward{\noexpand#1}{page}{\the\pageno}}%
   \bgroup\htmlanchor{#1}{\makeref{\noexpand#1}{page}{\the\pageno}}\egroup%
   \draftlabel{page\ #1}%
 }

\def\plot[#1]{
    \ifdraft\let\proclaimcolor=\Fuchsia\else\let\proclaimcolor=\Black\fi%
    \gdef\Prnm{Plot}%
    \global\advance\thmno by 1%
    \ifnumct{\bf \proclaimcolor{Plot\ \proctoks}}%
    \else\noindent{\bf \proclaimcolor{Plot}}%
    \fi%
  \htmlanchor{#1}{\makeref{\noexpand#1}{\Prnm}{\proctoks}}%
  \immediate\write\isauxout{\noexpand\forward{\noexpand#1}{\Prnm}{\proctoks}}%
  \draftlabel{#1}%
  }

\newcount\captionct \captionct=0
\outer\def\caption#1{{\global\advance\captionct by 1%
  \relax%
  \gdef\Irnm{Figure}%
  \edef\numtoks{\ifchap\the\chapno.\fi\ifsect\the\sectno.\fi\the\captionct}%
  \vskip0.5cm\vtop{\noindent{\bf Figure \numtoks.} #1}\vskip0.5cm%
  \futurelet\testchar\maybeoptionproclaim}}
\outer\def\caption#1{{\global\advance\thmno by 1%
  \relax%
  \gdef\Irnm{Figure}%
  \edef\proctoks{\ifchap\the\chapno.\fi\ifsect\the\sectno.\fi\the\thmno}%
  \vskip0.2cm\vtop{\narrower\narrower\noindent{\bf Figure \proctoks.} #1}\vskip0.5cm%
  \futurelet\testchar\maybeoptionproclaim}}

\def\eqlabel#1{
    \ifdraft\let\labelcolor=\Red\else\let\labelcolor=\Black\fi%
    \global\advance\thmno by 1%
    \edef\proctoks{\ifchap\number\chapno.\fi\ifsect\number\sectno.\fi\number\thmno}%
    \outer\expandafter\def\csname#1\endcsname{\proctoks}%
    \htmlanchor{#1}{\makeref{\noexpand#1}{}{\proctoks}}%
    \immediate\write\isauxout{\noexpand\forward{\noexpand#1}{Equation}{(\proctoks)}}%
    \quad\qquad\hfill\eqno\hbox{\hfill\rm\noexpand\labelcolor{(\proctoks)}\draftcmt{\labelcolor{#1}}}%
}
\def\iqlabel#1{
    \ifdraft\let\labelcolor=\Red\else\let\labelcolor=\Black\fi%
    \global\advance\thmno by 1%
    \edef\proctoks{\ifchap\number\chapno.\fi\ifsect\number\sectno.\fi\number\thmno}%
    \outer\expandafter\def\csname#1\endcsname{\proctoks}%
    \htmlanchor{#1}{\makeref{\noexpand#1}{}{\proctoks}}%
    \immediate\write\isauxout{\noexpand\forward{\noexpand#1}{Inequality}{(\proctoks)}}%
    \quad\qquad\hfill\eqno\hbox{\hfill\rm\noexpand\labelcolor{(\proctoks)}\draftcmt{\labelcolor{#1}}}%
}

\def\eqalabel#1{
    \ifdraft\global\let\labelcolor=\Red\else\global\let\labelcolor=\Black\fi%
    \global\advance\thmno by 1%
    \gdef\proctoks{\ifchap\number\chapno.\fi\ifsect\number\sectno.\fi\number\thmno}%
    \outer\expandafter\def\csname#1\endcsname{\proctoks}%
    \htmlanchor{#1}{\makeref{\noexpand#1}{}{\proctoks}}%
    \immediate\write\isauxout{\noexpand\forward{\noexpand#1}{equation}{(\proctoks)}}%
    \quad\qquad&\hfill\hbox{\hfill\rm\noexpand\labelcolor{(\proctoks)}\draftcmt{\labelcolor{#1}}}%
}

\def\label#1{\optionproclaim[#1]}
\def\label#1{%
    \ifdraft\let\labelcolor=\Red\else\let\labelcolor=\Black\fi%
    \edef\proctoks{\ifchap\the\chapno.\fi\ifsect\the\sectno.\fi\the\thmno}%
    \outer\expandafter\def\csname#1\endcsname{\proctoks}%
    \bgroup\htmlanchor{#1}{\makeref{\noexpand#1}{}{\proctoks}}\egroup%
    \immediate\write\isauxout{\noexpand\forward{\noexpand#1}{Figure}{\proctoks}}%
}

\proclamation{definition}{Definition}{1}
\proclamation{defin}{Definition}{1}
\proclamation{defins}{Definitions}{1}
\proclamation{lemma}{Lemma}{1}
\proclamation{lemmap}{Lemma}{2}
\proclamation{lemmapp}{Lemma}{3}
\proclamation{proposition}{Proposition}{1}
\proclamation{prop}{Proposition}{1}
\proclamation{theorem}{Theorem}{1}
\proclamation{thm}{Theorem}{1}
\proclamation{corollary}{Corollary}{1}
\proclamation{cor}{Corollary}{1}
\proclamation{conjecture}{Conjecture}{1}
\proclamation{proc}{Procedure}{1}
\proclamation{assumption}{Assumption}{1}
\proclamation{axiom}{Axiom}{1}
\proclamation{example}{Example}{0}
\proclamation{examples}{Examples}{0}
\proclamation{examplesnotation}{Examples and notation}{0}
\proclamation{exnoex}{Examples and non-examples}{0}
\proclamation{remark}{Remark}{0}
\proclamation{remarks}{Remarks}{0}
\proclamation{remarksexamples}{Remarks and examples}{0}
\proclamation{facts}{Facts}{0}
\proclamation{fact}{Fact}{0}
\proclamation{irrelevant}{Irrelevant}{0}
\proclamation{question}{Question}{0}
\proclamation{construction}{Construction}{0}
\proclamation{algorithm}{Algorithm}{0}
\proclamation{problem}{Problem}{0}
\proclamation{table}{Table}{0}
\proclamation{notation}{Notation}{0}
\proclamation{figure}{Figure}{0}
\proclamation{impexample}{Example}{-1}
\proclamation{exercise}{Exercise}{0}
\proclamation{importantexercise}{IMPORTANT EXERCISE}{0}
\proclamation{exercisedag}{Exercise}{7} 
\proclamation{stexercise}{Exercise}{8} 
\proclamation{excise}{Exercise}{6} 
\proclamation{stexcise}{Exercise}{9} 
\proclamation{excisedag}{Exercise}{10} 

\def\doubleindent{\advance \leftskip 2\parindent\advance \rightskip \parindent}
\def\singleindent{\advance \leftskip \parindent\advance \rightskip \parindent}
\def\doubleleftindent{\advance \leftskip 2\parindent}


\def\hpad#1#2#3{\hbox{\kern #1\hbox{#3}\kern #2}}
\def\vpad#1#2#3{\setbox0=\hbox{#3}\vbox{\kern #1\box0\kern #2}}



\def\stack#1#2#3{\vbox{\offinterlineskip
  \setbox2=\hbox{#2}
  \setbox3=\hbox{#3}
  \dimen0=\ifdim\wd2>\wd3\wd2\else\wd3\fi
  \hbox to \dimen0{\hss\box2\hss}
  \kern #1
  \hbox to \dimen0{\hss\box3\hss}}}


\def\hexp#1{%
  \setbox0=\hbox{${}^{#1}$}%
  \hbox to .5\wd0{\box0\hss}}

\def\hsub#1{%
  \setbox0=\hbox{${}_{#1}$}%
  \hbox to .5\wd0{\box0\hss}}


\def\mOth{\mathsurround=0pt}
\newdimen\pOrenwd \setbox0=\hbox{\tenex B} \pOrenwd=\wd0
\def\leftextramatrix#1{\begingroup \mOth
  \setbox0=\vbox{\def\cr{\crcr\noalign{\kern2pt\global\let\cr=\endline}}
    \ialign{$##$\hfil\kern2pt\kern\pOrenwd&\thinspace\hfil$##$\hfil
      &&\quad\hfil$##$\hfil\crcr
      \omit\strut\hfil\crcr\noalign{\kern-\baselineskip}
      #1\crcr\omit\strut\cr}}
  \setbox2=\vbox{\unvcopy0 \global\setbox1=\lastbox}
  \setbox2=\hbox{\unhbox1 \unskip \global\setbox1=\lastbox}
  \setbox2=\hbox{$\kern\wd1\kern-\pOrenwd \left[ \kern-\wd1
    \vcenter{\unvbox0 \kern-\baselineskip} \,\right]$}
  \null\;\vbox{\box2}\endgroup}


\def\COMMENT#1\par{\bigskip\hrule\smallskip#1\smallskip\hrule\bigskip}

{\catcode`\^^M=12 \endlinechar=-1 %
 \gdef\xcomment#1^^M{\def\test{#1}
   \ifx\test\endcomment \let\next=\endgroup
   \else \let\next=\xcomment \fi
   \next}
}
\def\dospecials{\do\ \do\\\do\{\do\}\do\$\do\&%
  \do\#\do\^\do\^^K\do\_\do\^^A\do\%\do\~}
\def\makeinnocent#1{\catcode`#1=12}
\def\comment{\begingroup
  \let\do=\makeinnocent \dospecials
  \endlinechar`\^^M \catcode`\^^M=12 \xcomment}
{\escapechar=-1
 \xdef\endcomment{\string\\endcomment}
}



\def\draftcmt#1{\ifdraft\llap{\smash{\raise2ex\hbox{{#1}}}}\fi}
\def\draftind#1{\llap{\smash{\raise-1.3ex\hbox{{#1}}}}}
\def\draftlabel#1{\ifdraft\draftcmt{\footfonts\Green{{\string#1}\ }\unskip\fi}}
\def\draftlabel#1{%
  \ifx#1\empty\else
  \ifdraft\draftcmt{\foott\Green{\string #1}\ \ }\unskip\fi\fi}

\newbox\strikebox
\def\strike#1{\setbox\strikebox=\hbox{#1}%
  \rlap{\raise0.4ex\hbox to \wd\strikebox{\leaders\hrule height 0.8pt depth 0pt\hfill}}%
  #1}
\def\strikem#1{\setbox\strikebox=\hbox{$#1$}%
  \rlap{\raise0.4ex\hbox to \wd\strikebox{\leaders\hrule height 0.8pt depth 0pt\hfill}}%
  \hbox{$#1$}}

\newcount\hours
\newcount\minutes
\hours\time \divide\hours 60
\minutes-\hours \multiply\minutes 60\advance\minutes\time
\edef\timehhmm{\ifnum\hours<10 0\fi\the\hours
:\ifnum\minutes<10 0\fi\the\minutes}

\footline={{\rm\hfill\folio\hfill}}

\def\cite{\futurelet\testchar\maybeoptioncite}
\def\maybeoptioncite{\ifx[\testchar \let\next\optioncite
	\else \let\next\nooptioncite \fi \next}

\def\nooptioncite#1{\expandafter\ifx\csname#1\endcsname\relax
  {\nodefmsg{#1}\bf`#1'}\else
  \htmllocref{#1}{[\csname#1\endcsname]}{}\fi}

\def\optioncite[#1]#2{\ifx\csname#1\endcsname\relax\else
  \htmllocref{#2}{[\csname#2\endcsname, #1]}{}\fi}


\def\strutdepth{\dp\strutbox}
\def\strutdepth{4.3ex}
\def\margin#1{{\strut\vadjust{\kern-\strutdepth%
        \vtop to\strutdepth{\baselineskip\strutdepth%
        \llap{\vbox{\hsize=3em\noindent #1}\hskip1em}\null}}}}
\def\rmargin#1{\strut\vadjust{\kern-\strutdepth%
        \vtop to\strutdepth{\baselineskip\strutdepth%
        \vss\rlap{\ \hskip\hsize\ \vtop{\hsize=4em\noindent #1}}\null}}}
\def\warning{\vskip0ex\noindent\llap{\lower4pt\hbox to 1.0cm{\psline(0,.3)(.27,.3)(0,.09)(.03,0)(.27,0)(.3,.03)(.27,.06)}\hfill\ }\indent}
\def\warning{\margin{\lower4pt\hbox to 0.6cm{\psline(0,.3)(.27,.3)(0,.09)(.03,0)(.27,0)(.3,.03)(.27,.06)}\hfill\ }}
\def\warning{\margin{\hbox to 0.6cm{\psline(0,.3)(.27,.3)(0,.09)(.03,0)(.27,0)(.3,.03)(.27,.06)}\hfill\ }}
\def\warning{\noindent{\bf Warning: }\margin{\ \ \ \ \lower-4pt\hbox to 0.6cm{\psline[linewidth=2pt](0,0)(.5,0)(0.53,-0.02)(0.55,-0.05)(0.53,-0.1)(0.03,-0.35)(0.01,-0.4)(0.03,-0.43)(0.06,-0.45)(0.56,-0.45)}\hfill}}
\def\warningdown#1{\noindent{\bf Warning: }\bgroup\def\strutdepth{#1}\margin{\ \ \
\ \hbox to 0.6cm{\psline[linewidth=2pt](0,0)(.5,0)(0.53,-0.02)(0.55,-0.05)(0.53,-0.1)(0.03,-0.35)(0.01,-0.4)(0.03,-0.43)(0.06,-0.45)(0.56,-0.45)}\hfill}\egroup}

\def\quote{\bgroup\footfont\parindent=0pt%
	\baselineskip=8pt\hangindent=15em\hangafter=0}
\def\endquote{\egroup}

\def\bmargin#1{{\strut\vadjust{\kern-\strutdepth%
        \vtop to\strutdepth{\baselineskip\strutdepth%
        \llap{\vbox{\hsize=3em\noindent #1}\hskip2em}\null}}}}

\newdimen\hngind \hngind=\parindent \advance\hngind by 3.9em\relax
\def\nonuritem{{\vskip0pt\hangafter=0\global\hangindent\hngind%
        {\vphantom{(\the\randomct)\hskip.5em}}}}
\def\nonuritem{{\parindent=0pt\hngind=\parindent \advance\hngind by 1.9em%
        \vskip0pt\hangafter=0\global\hangindent\hngind%
        \vphantom{{(\the\randomct)\hskip.5em}}}}
\def\nonuritem{{\idoitem{}}}
\def\pitem{{\global\advance\randomct by 1%
	\vskip1pt\hangafter=0\global\hangindent\hngind%
	{{\the\randomct.\hskip.5em}}}}

%
\def\idoitem#1{{\parindent=0pt\hngind=\parindent \advance\hngind by 3.7em%
        \vskip0pt\hangafter=0\global\hangindent\hngind%
        {$\hbox to3.7em{\hfill#1\hskip.5em}$}}}%
\def\romanitem{\global\advance\randomct by 1%
    \edef\itoks{(\romannumeral\randomct)}%
    \edef\Irnm{}%
    \idoitem{\itoks}%
    \futurelet\testchar\maybeoptioniitem} 
\def\sitem{\global\advance\randomct by 1%
    \edef\itoks{(\the\randomct)}%
    \edef\Irnm{Step}%
    \idoitem{\itoks}%
    \futurelet\testchar\maybeoptioniitem} 
\def\Rsitem{\global\advance\rrandomct by 1%
    \edef\itoks{(\uppercase\expandafter{\romannumeral\rrandomct)}}%
    \edef\Irnm{Step}%
    \idoitem{\itoks}%
    \futurelet\testchar\maybeoptioniitem} 
\def\citem{\global\advance\randomct by 1%
    \edef\itoks{(\the\randomct)}%
    \edef\Irnm{Condition}%
    \idoitem{\itoks}%
    \futurelet\testchar\maybeoptionitem} 
\def\iitem{\global\advance\randomct by 1%
    \edef\itoks{(\the\randomct)}%
    \edef\Irnm{}%
    \idoitem{\itoks}%
    \futurelet\testchar\maybeoptioniitem} 
\def\rpitem{\global\advance\randomct by 1%
    \edef\itoks{(\the\randomct')}%
    \edef\Irnm{}%
    \idoitem{\itoks}%
    \futurelet\testchar\maybeoptioniitem} 

\def\maybeoptioniitem{\ifx[\testchar\let\next\optioniitem%
	\else\let\next\nooptioniitem\fi\next}
\def\optioniitem[#1]{%
  \bgroup\htmlanchor{#1}{\makeref{\noexpand#1}{\Irnm}{\itoks}}\egroup%
  \immediate\write\isauxout{\noexpand\forward{\noexpand#1}{\Irnm}{\itoks}}%
  \draftlabel{#1}}
\def\nooptioniitem{}
\def\bulletpoint{{\parindent=0pt\hngind=\parindent \advance\hngind by 1.2em%
        \vskip0pt\hangafter=0\global\hangindent\hngind%
        {$\hbox to1.2em{\hfill$\bullet$\hskip.5em}$}}}%
\def\bulletpoint{{\hngind=\parindent \advance\parindent by -1.2em%
        \vskip0pt\hangafter=0\global\hangindent\hngind%
        {$\hbox to1.2em{\hfill$\bullet$\hskip.5em}$}}}%
\def\nopoint{{\parindent=0pt\hngind=1cm \advance\hngind by 1.2em%
        \vskip0pt\hangafter=0\global\hangindent\hngind%
        {$\hbox to1.2em{\hfill\ \hskip.5em}$}}}%
\def\nopoint{{\hngind=\parindent \advance\parindent by -1.2em%
        \vskip0pt\hangafter=0\global\hangindent\hngind%
        {$\hbox to1.2em{\hfill\ \hskip.5em}$}}}%
\def\bbulletpoint{{\bgroup\parskip=0.1cm\parindent=1.5cm\hngind=\parindent%
	\advance\hngind by 1.5em%
        \vskip0pt\hangafter=0\global\hangindent\hngind%
        {$\hbox to1.5em{\hfill$\bullet$\hskip.5em}$}\egroup}}%
\def\bnobulletpoint{{\bgroup\parskip=0.1cm\parindent=1.5cm\hngind=\parindent%
	\advance\hngind by 1.5em%
        \vskip0pt\hangafter=0\global\hangindent\hngind%
        {$\hbox to1.5em{\hfill\hskip.5em}$}\egroup}}%

\let\ritem=\iitem

\def\rlabel#1{{%
  \bgroup\htmlanchor{#1}{\makeref{\noexpand#1}{item}{\the\randomct}}\egroup%
  \immediate\write\isauxout{\noexpand\forward{\noexpand#1}{item}{\the\randomct}}%
  \draftlabel{item\ #1}}}

\def\idoitemitem#1{{\hngind=\parindent \advance\hngind by 4.3em%
        \vskip-\parskip\hangafter=0\global\hangindent\hngind%
        {$\hbox to3.9em{\hfill#1\ }$}}}
\def\idoitemitemm#1{{\hngind=\parindent \advance\hngind by 7.3em%
        \vskip-\parskip\hangafter=0\global\hangindent\hngind%
        {$\hbox to6.9em{\hfill#1\ }$}}}


\def\thmparens#1{\ignorespaces{\rm(#1)}}

\newdimen\thmskip \thmskip=0.7ex
	
\def\thmalready#1{\removelastskip\vskip\thmskip%
	\global\randomct=0%
	\noindent{\bf #1:}\ \ %
	\bgroup \it%
	\abovedisplayskip=4pt\belowdisplayskip=3pt%
	\parskip=0pt%
	}


\newdimen\exerciseindent \exerciseindent=4.5em
\newdimen\exerciseitemindent \exerciseitemindent=1.9em
\newdimen\exitindwk \exitindwk=3.7em
\def\exitem{\global\advance\randomct by 1%
	\vskip0pt\indent\hangindent\exitindwk \hskip\exerciseitemindent%
	\llap{\romannumeral\randomct)\hskip0.5em}\ \hskip-0.5em}
\def\stexitem{\global\advance\randomct by 1%
	\par\indent\hangindent\exitindwk \hskip\exerciseitemindent%
	\llap{(\romannumeral\randomct)\hskip-.1ex*\enspace}\ignorespaces}


\def\exercises{\global\thmno=0\vskip0.9cm%
	\noindent{\bf Exercises for Section~\ifchap\the\chapno.\fi\the\sectno}}
\def\solutions{\global\thmno=0\bigskip\widow{0.05}%
	\noindent{\bf Solutions for Section~\ifchap\the\chapno.\fi\the\sectno}}
\def\solution#1{\removelastskip\vskip 5pt%
	\global\randomct=0%
	\parindent=\exerciseindent%
	\par\hangindent\exerciseindent\indent\llap{\hbox to \exerciseindent{%
	  {\bf \unskip\ref{#1}}}}}
\def\solution#1{\removelastskip\vskip 5pt\par\noindent {\bf \refn{#1}:}}


\newdimen\wddim
\def\edinsert#1{\setbox0=\hbox{#1}\wddim=\wd0\divide\wddim by 2%
\smash{\hskip-\wddim\raise15pt\hbox to 0pt{\Red{$\underbrace{\box0}$}}\hskip\wddim}}

\def\midline#1{\raise 0.20em\hbox to #1{\vrule depth0pt height 0.405pt width #1}}

\def\buildrelu#1\over#2{\mathrel{\mathop{\kern0pt #1}\limits_{#2}}}



\def\bmatrix#1#2{{\left(\vcenter{\halign
  {&\kern#1\hfil$##\mathstrut$\kern#1\cr#2}}\right)}}

\def\rightarrowmat#1#2#3{
  \setbox1=\hbox{\small\kern#2$\bmatrix{#1}{#3}$\kern#2}
  \,\vbox{\offinterlineskip\hbox to\wd1{\hfil\copy1\hfil}
    \kern 3pt\hbox to\wd1{\rightarrowfill}}\,}

\def\leftarrowmat#1#2#3{
  \setbox1=\hbox{\small\kern#2$\bmatrix{#1}{#3}$\kern#2}
  \,\vbox{\offinterlineskip\hbox to\wd1{\hfil\copy1\hfil}
    \kern 3pt\hbox to\wd1{\leftarrowfill}}\,}

\def\rightarrowbox#1#2{
  \setbox1=\hbox{\kern#1\hbox{\small #2}\kern#1}
  \,\vbox{\offinterlineskip\hbox to\wd1{\hfil\copy1\hfil}
    \kern 3pt\hbox to\wd1{\rightarrowfill}}\,}

\def\leftarrowbox#1#2{
  \setbox1=\hbox{\kern#1\hbox{\small #2}\kern#1}
  \,\vbox{\offinterlineskip\hbox to\wd1{\hfil\copy1\hfil}
    \kern 3pt\hbox to\wd1{\leftarrowfill}}\,}


\newcount\countA
\newcount\countB
\newcount\countC

\def\monthname{\begingroup
  \ifcase\number\month
    \or January\or February\or March\or April\or May\or June\or
    July\or August\or September\or October\or November\or December\fi
\endgroup}

\def\dayname{\begingroup
  \countA=\number\day
  \countB=\number\year
  \advance\countA by 0 
  \advance\countA by \ifcase\month\or
    0\or 31\or 59\or 90\or 120\or 151\or
    181\or 212\or 243\or 273\or 304\or 334\fi
  \advance\countB by -1995
  \multiply\countB by 365
  \advance\countA by \countB
  \countB=\countA
  \divide\countB by 7
  \multiply\countB by 7
  \advance\countA by -\countB
  \advance\countA by 1
  \ifcase\countA\or Sunday\or Monday\or Tuesday\or Wednesday\or
    Thursday\or Friday\or Saturday\fi
\endgroup}

\def\timename{\begingroup
   \countA = \time
   \divide\countA by 60
   \countB = \countA
   \countC = \time
   \multiply\countA by 60
   \advance\countC by -\countA
   \ifnum\countC<10\toks1={0}\else\toks1={}\fi
   \ifnum\countB<12 \toks0={\sevenrm AM}
     \else\toks0={\sevenrm PM}\advance\countB by -12\fi
   \relax\ifnum\countB=0\countB=12\fi
   \hbox{\the\countB:\the\toks1 \the\countC \thinspace \the\toks0}
\endgroup}

\def\timestamp{\dayname, \the\day\ \monthname\ \the\year, \timename}

%

%


\def\frac#1#2{{#1 \over #2}}

\let\text\hbox
\def\overset#1\to#2{\ \mathop{\buildrel #1 \over #2}\ }
\def\underset#1\to#2{{\mathop{\buildrel #1 \over #2}}}

\def\myoperator{NOOP}
\def\mymathopsp{\futurelet\testchar\maybesubscriptop}
\def\maybesubscriptop{\ifx_\testchar \let\next\subscriptop
	\else \let\next\nosubscriptop \fi \next}
\def\subscriptop_#1{{\rm \myoperator}_{#1}}
\def\nosubscriptop{\futurelet\testchar\maybeparenop}
\def\maybeparenop{{\mathop{\rm \myoperator}\nolimits%
	\ifx(\testchar \relax\else\,\fi}}


\def\notdiv{\hbox{$\not\hbox to -0.2ex{$|$}$}\hskip0.5em}

\def\Ass{\mathop{\rm Ass}\nolimits}



\def\ker{\mathop{\rm ker}\nolimits}

\def\semidirect{{\mathop{\hbox to 1.1em{\hskip.5em\psline[linewidth=0.3pt](-.18,+.01)(.18,.15)(.18,.01)(-.18,.15)}}}}


\def\spm{\raise0.5pt\hbox{\bgroup\textfont2=\fivesy\relax$\pm$\egroup}}

\def\bbZ{{\hbox{\tenbbold Z}}}

\def\bbN{{\hbox{\tenbbold N}}}


\def\uspace{\hbox{$\underline{\hbox to 1ex{\ \hfil}}$}}

\def\lineunderblank{\underline{\ \hbox to 0.5ex{$\vphantom{a}$\hfill\ }}}
\def\lineunderblankb#1{\underline{\ \hbox to #1{$\vphantom{a}$\hfill\ }}}
\def\lineaboveblank{\overline{\ \hbox to 0.5ex{$\vphantom{a}$\hfill\ }}}
\def\lineaboveblank{\overline{\ \hbox to 0.5ex{$\vphantom{a}$\hfill\ }}}
\def\blank#1{\hbox{$\underline{\hbox to #1em{\ \hfill}}$}}


\def\qedbox{\hbox{\vbox{\hrule\hbox{\vrule\kern3pt\vbox{\kern6pt}\kern3pt\vrule}\hrule}}}
\def\qed{\unskip\hfill\qedbox\vskip\thmskip} 
\def\eqed{\eqno\hbox{\quad\qedbox}} 


\def\blankpage{\vfill\eject\bgroup\footline={}\headline={}\ \vtop to0cm{
\pspolygon[linecolor=white](0,0)(0,1)(1,0)\vfill}
\vfill\eject\egroup}
\def\blankpage{\vfill\eject\bgroup\footline={}\headline={}\ \vtop to0cm{
\ \vfill}\vfill\eject\egroup}


\def\xspace{\futurelet\lettoken\doxspace}
\def\doxspace{%
  \ifx\lettoken\bgroup\else%
  \ifx\lettoken\egroup\else%
  \ifx\lettoken\/\else%
  \ifx\lettoken\ \else%
  \ifx\lettoken~\else%
  \ifx\lettoken.\else%
  \ifx\lettoken!\else%
  \ifx\lettoken,\else%
  \ifx\lettoken:\else%
  \ifx\lettoken;\else%
  \ifx\lettoken?\else%
  \ifx\lettoken/\else%
  \ifx\lettoken'\else%
  \ifx\lettoken)\else%
  \ifx\lettoken-\else%
  \ifx\lettoken\space\else%
   \space%
   \fi\fi\fi\fi\fi\fi\fi\fi\fi\fi\fi\fi\fi\fi\fi\fi}


\def\today{\ifcase\month\or January\or February\or March\or
  April\or May\or June\or July\or August\or September\or
  October\or November\or December\fi
  \space\number\day, \number\year}\relax


\newwrite\isauxout
\openin1\jobname.aux
\ifeof1\message{No file \jobname.aux}
       \else\closein1\relax\input\jobname.aux
       \fi
\immediate\openout\isauxout=\jobname.aux
\immediate\write\isauxout{\relax}

\newwrite\iscontout
\immediate\openout\iscontout=\jobname.cont
\immediate\write\iscontout{\relax}
\def\chaplist#1#2#3{{\bf \vskip1pt\noindent\ifx0#1#2\hfill#3%
		\else\ref{#1}: #2\hfill#3%
                \fi}}
\def\basicchaplist#1#2#3{{\bf \vskip1pt\noindent\ifx0#1#2\hfill#3%
		\else#2\hfill#3%
                \fi}}
\def\contlist#1#2#3{\ifchap\vskip1pt\noindent \ \ \ifx0#1#2\hfill#3%
		\else\ref{#1}: #2\hfill#3%
                \fi%
	\else \vskip1pt\noindent\ref{#1}: #2\hfill#3\fi}
\def\contlist#1#2#3{\vskip1pt\noindent \ \hskip1em%
	\ifchap\ifx0#1#2\hfill#3%
		\else\ref{#1}: #2\hfill#3%
                \fi%
	\else \ref{#1}: #2\hfill#3\fi}
\def\contnosectlist#1#2#3{\vskip1pt\noindent\hskip1em#2\hfill#3}


%
\newwrite\isindexout
\immediate\openout\isindexout=\jobname.index
\font\footttt=cmtt6

\def\padno{
\ifnum\pageno<10 00\the\pageno\else\ifnum\pageno<100 0\the\pageno\else\the\pageno\fi\fi}

\def\index#1{\rlap{\ifdraft\draftind{{\footttt\Blue{#1\ \ }}}\unskip\fi%
\write\isindexout{\noexpand\indexentry {#1}{\padno}}}\unskip\ignorespaces}

\def\indexn#1{\rlap{\ifdraft\draftind{{\footttt\Brown{#1\ \ }}}\unskip\fi%
\write\isindexout{\noexpand\indexentryn {#1}{\padno}}}\unskip\ignorespaces}

\def\obeyspaces{\catcode`\ =\active\catcode`\	=\active}
{\obeyspaces\global\let =\space} 
{\obeyspaces\gdef {\hskip.5em}\gdef	{\hskip3em}}
{\catcode`\^^M=\active %
\gdef\obeylines{\catcode`\^^M=\active \gdef^^M{\vskip.0ex\ }}}%

\font\stt=cmtt10
\def\dMaccode{\vskip3ex\bgroup\advance\leftskip1ex\parindent=0pt%
	\baselineskip=9pt\obeylines\obeyspaces%
	\catcode`\^=11\catcode`\_=11\catcode`\#=11\catcode`\%=11%
	\catcode`\{=11\catcode`\}=11%
	\stt}
\def\enddMaccode{\removelastskip\unskip\egroup}
\def\Maccode{\bgroup \catcode`\^=11 \catcode`\_=11 \stt"}
\def\endMaccode{\unskip"\egroup}
\def\Maccode{\bgroup \catcode`\^=11 \catcode`\_=11 \stt}
\def\endMaccode{\unskip\egroup}

\gdef\poetry{\vskip0pt\bgroup\parindent=0pt\obeylines\obeyspaces}
\gdef\poetrycenter{\vskip0pt\bgroup\parskip=0pt\leftskip=0pt plus 1fil\rightskip=0pt plus
1fil\parfillskip=0pt\obeylines}
\gdef\endpoetry{\vskip0pt\egroup}

\irenasize

\newwrite\isauxout
\openin1\jobname.aux
\ifeof1\message{No file \jobname.aux}
       \else\closein1\relax\input\jobname.aux
       \fi
\immediate\openout\isauxout=\jobname.aux
\immediate\write\isauxout{\relax}

\def\spread{\sigma}
\def\spread{spr}
\def\spread{\mathop{\rm spr}\nolimits}

\drafttrue
\draftfalse
\colortrue
\ifcolor\input colordvi\else\input blackdvi\fi%
\ifdraft\footline={{\rm \hfill page \folio -- \today\ --
\timehhmm\hfill}}\else
\ifnum\pageno>1\footline={{\rm\hfill\folio\hfill}}\else\footline={}\fi\fi

\centerline{\sectionfont Many associated primes of powers of primes}
\centerline{\sectionfont Jesse Kim and Irena Swanson}
\centerline{Department of Mathematics, UCSD, La Jolla, CA 92093-0112}
\centerline{Reed College, 3203 SE Woodstock Boulevard, Portland, OR 97202}
\centerline{\tt jvkim@ucsd.edu, iswanson@reed.edu}

\vskip 0.5cm
\bgroup
\narrower\narrower
\noindent
{\bf Abstract.}
\sl
We construct families of prime ideals in polynomial rings
for which
the number of associated primes of the second power (or higher powers)
is exponential in the number of variables in the ring.

\vskip 0.3cm
\noindent
{\it Keywords and phrases:}
primary decomposition,
powers of ideals,
Rees algebra,
extended Rees algebra,
Rees-like algebra,
polynomial rings.

\noindent
{\it 2010 Mathematics Subject Classification. Primary 13A30, 13B25.}


\egroup 
\vskip 0.5cm

An inspiration for this work came from McCullough-Peeva's paper~\cite{McP}
in which they constructed
families of counterexamples to the Eisenbud-Goto conjecture
via new constructions
of step-by-step homogenizations and Rees-like algebras.
McCullough and Peeva applied these constructions in particular to the
Mayr-Meyer ideals from~\cite{MM}.
Mayr-Meyer ideals are computationally hard
in the sense that they have large Castelnuovo-Mumford regularity
(Bayer and Stillman~\cite{BS})
and large ideal-membership coefficient degrees (Mayr and Meyer~\cite{MM}).
Another class of computationally hard ideals are permanental ideals
(\cite{V}).
Papers~\cite{K}, \cite{LS}, \cite{PS},
\cite{S2}, \cite{S3}
indicate that
perhaps computational hardness
is related to the large numbers of associated primes,
but McCullough and Peeva showed that even prime ideals
--- so ideals with only one associated prime ---
can be computationally hard,
namely that they can have very large Castelnuovo-Mumford regularity.
This paper is a result of trying to understand
why McCullough and Peeva's
Rees-like algebras and step-by-step homogenizations
generate such ``hard" prime ideals,
and trying to determine whether
the primes defining Rees-like algebras
have large numbers of associated primes.
While initially we worked with Rees-like algebras,
we got tighter results with extended Rees algebras.

We construct classes of prime ideals~$P$ in polynomial rings over fields
such that the number of associated primes of $P^2$
is not bounded by any polynomial function in the number of variables.
In \ref{thmmainone}
we construct a family of
almost complete intersection prime ideals $P$ of height $n-1$
in $3n$ variables
for which all $P^e$ with $e \ge 2$
have exactly $3^n$ embedded primes.
These $P$ are generated by elements of degrees up to $3$
(or by quasi-homogeneous elements of degrees up to $10$).
In \ref{thmmainmain} we construct for odd $n$
a family of prime ideals $P$ of height ${n+1 \over 2}$ in $3n$ variables
for which $P^2$ has at least
$3^n + 3^{(n+7)/2}$ embedded primes.
We know of no other class of ideals
with such large numbers of associated prime ideals.
The second author proved in~\cite{S3}
that there is an upper bound on the number
of associated prime ideals of the Mayr-Meyer ideals
that is doubly exponential in the number of variables,
but it is not known whether the number of associated primes
is in fact doubly or even singly exponential.

Hermann \cite{H} was the first to consider upper bounds
on the numbers and degrees of primary components of
ideals in polynomial rings over a field.
Seidenberg proved in~\cite[Point 65]{Sei74}
that there exists a primitive recursive function $B(n,d)$
that is at least doubly exponential in $n$
such that 
any ideal~$I$ in a polynomial ring in $n$~variables over a field
with generators of degree at most~$d$
has at most $B(n,d)$ associated primes.
Our \ref{thmmainone} shows that $B(n,d) \ge~3^{n/3}$ for all $d \ge 3$.
More recent proofs of upper bounds on the numbers of primary components
are in the paper~\cite{vdDS} by van den Vries and Schmidt
and in the paper~\cite{AH} by Ananyan and Hochster.
The bound $E(m,d)$ given by Ananyan and Hochster depends on 
the upper bound $m$ on the generators of $I$
and on the upper bound $d$ of the degrees of the generators,
and it does not depend on the number of variables.
(Note that $B(n,d)$, while not explicitly invoking $m$,
does get a free upper bound ${n+d \choose d}$ on $m$.)
In \ref{thmmainone}
we construct a family of
$(n+1)$-generated almost complete intersection prime ideals~$P$
in $3n$ variables
for which all $P^e$ with $e \ge 2$
have exactly $3^n$ embedded primes.
Thus $P^2$ has at most $m = (n+2)(n+1)/2 \le (n+2)^2/2$ generators,
which shows that $E(m,d) \ge 3^{\sqrt{2m}-1}$.

We develop two new methods for generating ideals with
large numbers of primary components:
splitting, which is a generalization of the step-by-step homogenization,
and spreading.
\ref{sectsplflat} develops some basic properties of splitting.
We show that splitting increases the number of associated primes
in a controlled way;
the number is bounded above by a polynomial
in the number of variables in the new ring,
and the polynomial is of degree that is equal to the largest number
of variables contained in an associated prime of the original ideal.
Splitting does not increase the number of variables in the associated primes;
it increases the numbers of associated primes.
With the goal of increasing the numbers of variables in associated primes
we introduce in \ref{sectspreading} the new notion of spreading.
We determine the presenting ideal
of the Rees algebra, of the extended Rees algebra, or of the Rees-like algebra
of the spreading of an ideal $I$
from the presenting ideal of the same type of Rees algebra of~$I$.

In \ref{sectexamples} we compute the presenting ideal~$P$
of the extended Rees algebra of a specific five-generated monomial ideal.
The spreading of monomial ideals is easy to understand.
In \ref{sectsplexamples} we apply spreading and splitting to this~$P$
to get the exponential results.

All rings in this paper are commutative with identity,
and most are Noetherian.

\bigskip
{\bf Acknowledgements.}
We are grateful to the two anonymous referees for substantially
improving the presentation and the proofs,
in particular of \ref{thmsplittingff}
and of \refs{lmprimprimary} and \refn{lmsplitprimwitht}.
The paper is now greatly streamlined due to their comments,
and the main results are stated without reliance on calculations by Macaulay2.

\section{Splitting, flatness, primary decompositions}[sectsplflat]%

In this section we define splittings,
we prove that they are faithfully flat maps,
and we show that under splitting
the number of associated primes increases.
We give a lower bound for the number of associated primes,
and we show that the bound is achieved
when all exponents of variables in the splitting are equal to~$1$.

\defin
Let $A$ be a ring.
An {\bf $\bf A$-splitting}
is an $A$-algebra homomorphism
$\varphi : A[z] \to A[u_1, \ldots, u_n]$
given by $\varphi (z) = u_1^{p_1} \cdots u_n^{p_n}$,
where $n, p_1, \ldots, p_n$ are positive integers
and $z, u_1, \ldots, u_n$ are variables over $A$.
We refer to this map as the $A$-splitting
$z \mapsto u_1^{p_1} \cdots u_n^{p_n}$.
\endb

Step-by-step homogenization of McCullough-Peeva~\cite{McP}
is a special case of splitting
with $n = 2$, $p_1 = 1$ and $p_2$ chosen carefully depending on the gradings.
In this paper the splittings ignore any gradings.

The composition of splittings is a splitting.
If we restrict the splittings to those
for which at least one (resp.\ each) $p_i$ equals $1$,
then again the composition of two such splittings is of the same type.

\thm[thmsplittingff]
Splitting is a free and thus a faithfully flat map.
\endb

\proof
Note that the $A$-splitting $z \mapsto u_1^{p_1} \cdots u_n^{p_n}$
is a composition of the $A$-splitting
$A[z] \to A[u_1^{p_1}, \ldots, u_n^{p_n}]$
with the inclusion into the free extension $A[u_1, \ldots, u_n]$.
So it suffices to prove that the $A$-splitting
$A[z] \to A[u_1^{p_1}, \ldots, u_n^{p_n}]$ is free.
But $u_1^{p_1}, \ldots, u_n^{p_n}$ are variables over $A$,
so by possibly renaming them,
it remains to prove that the splitting $z \mapsto u_1 \cdots u_n$ is free.
But $A[u_1, \ldots, u_n]$ is free over $A[u_1 \cdots u_n] \cong A[z]$
with basis consisting of monomials $u_1^{i_1} u_2^{i_2} \cdots u_n^{i_n}$
with $i_1, i_2, \ldots, i_n$ non-negative integers
of which at least one equals $0$.
\qed

\lemma[lmprimprimary]%
\thmparens{This is a generalization of \cite[Exercise 10.4]{Eisenbud}.}
Let $a, b$ be a regular sequence in a ring $R$
and let $u$ be a variable over $R$.
Then $a$ is a non-zerodivisor on $R[u]/(au-b)R[u]$.
If the zero ideal in $R$ is prime (resp.\ primary),
then $(au-b)R[u]$ is a prime (resp.\ primary) ideal in $R[u]$.
\endb

\proof
Let $r \in (au-b) : a$.
Then $ra = s (au-b)$ for some $s \in R[u]$.
Since $a, b$ is a regular sequence in $R[u]$,
this means that there exists $s' \in R[u]$
such that $r-su = s' b$ and $s = - s' a$.
Hence $r = su + s'b = -s'(au - b) \in (au-b)$.
This proves the first statement.

Let $B = R[u]/(au-b)R[u]$.
Since $a$ is a non-zerodivisor on $B$ and on $R$,
we can form localizations $R_a$ and $B_a$
at the multiplicatively closed set $\{1,a,a^2,\ldots\}$.
Then $B$ injects into $B_a = R_a[u]/(u-b/a)R_a[u] \cong R_a$.
Since the zero ideal in $R$ is prime (resp.\ primary),
it is prime (resp.\ primary) in $R_a$
and hence also in $B_a$ and in its subring $B$.
This says that $(au-b)R[u]$ is a prime (resp.\ primary) ideal in $R[u]$.
\qed

\remark
It is not true in general with the set-up as in the lemma
that for an integer $p > 1$,
$(au^p - b)A[u]$ is prime or primary.
For example,
if $a = c^p, b = d^p$ is a regular sequence for some $c, d \in R$,
then $au^p - b$ factors.

\lemma[lmprimsplitnot]%
Let $A$ be a ring
and let $\varphi : A[z] \to A[u_1, \ldots, u_n]$
be the splitting map $z \mapsto u_1 \cdots u_n$.
Let $q$ be a prime (resp.\ primary) ideal in $A[z]$
such that $z$ is not in the radical of $q$.
Then $\varphi(q)A[u_1, \ldots, u_n]$ is prime (resp.\ primary)
of the same height as $q$,
and the $u_i$ are non-zerodivisor on $A[u_1, \ldots,
u_n]/\varphi(q)A[u_1,\ldots, u_n]$.
\endb

\proof
Let $U$ stand for $u_1, \ldots, u_n$.
Set $R = (A[z]/q)[U]$.
By assumption, $u_2\cdots u_n, z$ is a regular sequence on $R$.
Since $u_1$ is a variable in $R$ over the obvious subring,
we may apply \ref{lmprimprimary}
with $b = z, a = u_2 \cdots u_n, u = u_1$.
Then by lifting we get that
$qA[z][U] + (z-u_1\cdots u_n)$ is a prime (resp.\ primary) ideal in $A[z][U]$,
so that $\varphi(q)A[U] = (qA[z][U] + (z-u_1\cdots u_n)) \cap A[U]$
is a prime (resp.\ primary) ideal in $A[U]$.
\ref{lmprimprimary} also says that the $u_i$ are non-zerodivisors
on $R/(z-u_1\cdots u_n)$,
which by contraction means that they are non-zerodivisors
on $A[U]/\varphi(q)A[U]$.

The height claim follows from faithful flatness of splittings.
\qed

\lemma[lmsplitprimwitht]%
Let $A$ be a Noetherian ring
and let $\varphi : A[z] \to A[U] = A[u_1, \ldots, u_n]$
be the splitting map $z \mapsto u_1 \cdots u_n$.
Let $q$ be a primary ideal in $A[z]$ that contains a power of~$z$.
Then
$$
\varphi(q)A[U]
= \bigcap_{i=1}^n \biggl(\varphi(q)A[U] :
(u_1 \cdots u_{i-1} u_{i+1} \cdots u_n)^\infty\biggr)
$$
is an irredundant primary decomposition.
(Exponent $\infty$ stands for a very large integer;
the colon ideals are independent of the large integer
because the ring is Noetherian.)

The radical of $q$
can be written as $JA[z] + (z)$ with $J$ a prime ideal in $A$.
Then the associated primes of $\varphi(q)A[U]$
are $J A[U] + (u_1)$,
$\ldots$,
$J A[U] + (u_n)$.
In particular,
the heights of $q$ and $\varphi(q)A[U]$ are the same.
\endb

\proof
Let $Q$ be $\sqrt q$.
Then $Q$ is a prime ideal in $A[z]$ containing $z$,
so we can write it as $JA[z] + (z)$ for some prime ideal $J \subseteq A$.
Then $\varphi(Q)A[U] = JA[U] + (u_1\cdots u_n)$.
Set $Q_i = JA[U] + (u_i)$.
These are the prime ideals in $A[U]$
minimal over $\varphi(Q)A[U]$
and
$$
\bigcap_{i=1}^n Q_i = JA[U] + (u_1 \cdots u_n) = \varphi(Q)A[U]
$$
is an irredundant primary decomposition.
Since $\varphi$ is flat by \ref{thmsplittingff},
an application of \cite[Theorem 23.2]{Mats}
says that the associated primes of $qA[U]$ are precisely $Q_1, \ldots, Q_n$.
Since $u_j \in Q_j \setminus Q_i$ for all distinct $i, j$,
we get that the $Q_i$-primary component component of $\varphi(q)A[U]$
is $\varphi(q)A[U] : (u_1 \cdots u_{i-1} u_{i+1} \cdots u_n)^\infty$.
\qed

It should be noted that splitting
does not increase 
the numbers of variables in the associated primes.

\thm[thmpdgensplitmorevars]%
Let $A$ be a Noetherian ring,
$m, n_1, \ldots, n_m$ positive integers,
and let $A[Z] = A[z_1, \ldots, z_m]$
and $A[U] = A[u_{ij}: i = 1, \ldots, m; j = 1, \ldots, n_i]$
be polynomial rings over $A$.
Let $\varphi: A[Z] \to A[U]$
be the $A$-algebra homomorphism
with $z_i \mapsto u_{i1}^{p_{i1}} \cdots u_{in_i}^{p_{in_i}}$
for some positive integers $p_{ij}$.
Let $I$ be an ideal in $A[Z]$
with an irredundant primary decomposition
$I = q_1 \cap q_2 \cap \cdots \cap q_s$.
Then the following statements hold.
\ritem
The height of $\varphi(I) A[u_1, \ldots, u_n]$ equals the height of $I$.
\ritem
The number of primary components of $\varphi(I)A[U]$
is the sum of the numbers of primary components of the $\varphi(q_i)A[U]$.
\ritem
Let $I$ be a primary ideal in $A$.
For $i = 1, \ldots, m$
let $\epsilon_i$ be $1$ if $I$ contains a power of $z_i$
and let $\epsilon_i$ be $0$ otherwise.
Then
$\varphi(I)A[U]$ has at least $n_1^{\epsilon_1} \cdots n_m^{\epsilon_m}$
primary components.
If $p_{ij} = 1$ for all $(i,j)$ with $\epsilon_i = 1$,
then $\varphi(I)A[U]$ has exactly
$n_1^{\epsilon_1} \cdots n_m^{\epsilon_m}$ primary components.
The counted corresponding associated prime ideals
are of the form
$\sqrt {I} \cap A + (u_{ij} : \epsilon_i > 0, j = 1, \ldots, n_i)$.
\endb

\proof
The homomorphism $\varphi$
is a composition of the $A$-algebra homomorphisms
where each $z_i$ is split separately.
By \ref{thmsplittingff},
$\varphi$ is faithfully flat,
so that (1) holds.
Also,
flatness and \cite[Theorem 23.2]{Mats} say that
an irredundant primary decomposition of $\varphi(I) A[U]$
equals the intersection
of irredundant primary decompositions of the $\varphi(q_i) A[U]$.
This proves (2).

Suppose that (3) holds in case all $p_{ij} = 1$.
Then
$I' = I A[z_i, u_{ij}^{p_{ij}} : i, j]/
(z_i - u_{i1}^{p_{i1}} \cdots u_{in_i}^{p_{in_i}} : i)$
has the stated associated primes.
But $A[z_i, u_{ij}^{p_{ij}} : i, j]/
(z_i - u_{i1}^{p_{i1}} \cdots u_{in_i}^{p_{in_i}} : i)
\subseteq A[U] \cong A[z_i, u_{ij} : i, j]/
(z_i - u_{i1}^{p_{i1}} \cdots u_{in_i}^{p_{in_i}} : i)$
is a free extension,
so that
an irredundant primary decomposition of $\varphi(I) A[U]$
contracts to a possibly redundant primary decomposition of $I'$.
This means that the number of associated primes of $\varphi(I) A[U]$
is at least the number of primary components of $I'$.
Thus it suffices to prove (3) in case all $p_{1j}$ equal~$1$.
The case $m = 1$ follows from the previous two lemmas.

Now let $m > 1$.
Let $\varphi_1$ be the splitting of $z_1$
and let $\varphi'$ be the splitting of $z_2, \ldots, z_m$
such that $\varphi = \varphi' \circ \varphi_1$.
By the case $m = 1$,
the number of primary components of
$I' = \varphi_1(I)A[u_{11}, \ldots, u_{1n_1}]$ is $n_1^{\epsilon_1}$.
By the previous two lemmas,
for each $i > 1$ and for each primary component $q$ of $I'$,
$q$ contains $z_i$ if and only if $I$ contains $z_i$.
Thus by induction on~$m$,
the number of primary components of $\varphi'(q)A[U]$
is $n_2^{\epsilon_2} \cdots n_m^{\epsilon_m}$,
with the stated form.
By applying \cite[Theorem 23.2]{Mats} again due to faithful flatness,
$\varphi(I) A[U] = \varphi'(I')A[U]$ has
$n_1^{\epsilon_1} \cdots n_m^{\epsilon_m}$ primary components.
The corresponding associated primes are in the stated form.
\qed

\section{Variable spreading}[sectspreading]%

Splitting replaces one variable by a product of variables,
it can be described by a homomorphism,
and its effect is that after splitting an ideal,
the number of associated prime ideals increases.
By \ref{thmpdgensplitmorevars},
the increase is limited by the number of variables in the
associated primes of~$I$,
and it does not increase the number of variables in the associated primes.
This section introduces a new method,
which we call spreading,
which adds more variables.
A simple example is in \ref{propincreasevarseasy},
a more involved one in \ref{thmspreadReescommute}.

Here is an attempt at adding variables to associated primes.
Let $P$ (resp.~$Q$)
be a prime ideal in a polynomial ring $k[x_1, \ldots, x_n]$
(resp.~$k[y_1, \ldots, y_m]$)
with the property that it contains no variables but
all of its higher powers have embedded primes containing variables.
(Such an example is in \ref{exsimpleP}.)
Suppose that $(P+Q)k[x_1, \ldots, x_n, y_1, \ldots, y_m]$ is prime
(say if $k$ is algebraically closed).
Let $P^2 = P^{(2)} \cap p_1 \cap \cdots \cap p_r$
and
$Q^2 = Q^{(2)} \cap q_1 \cap \cdots \cap q_s$
be primary decompositions.
It turns out that the number of variables in the associated primes of
$(P+Q)^2$ does not increase
because $(P+Q)^{(2)}
\cap \bigcap_{i=1}^r (p_i + Q)
\cap \bigcap_{j=1}^s (P + q_j)$
is a primary decomposition of $(P+Q)^2$.
Namely,
by Theorem~2.7 in Walker~\cite{W},
$(P+Q)^{(2)} = P^{(2)} + PQ + Q^{(2)}$,
with the algebraic closure assumption
the ideals $p_i + Q$ and $P + q_j$ are primary,
and the intersection of these primary ideals is $(P+Q)^2$
by using that $\otimes_k k[\underline x]/P$
and $\otimes_k k[\underline y]/Q$ are flat.
Thus the associated primes of $(P+Q)^2$
are $P+Q$, $\sqrt {p_i} + Q$, $P + \sqrt{q_j}$,
and thus
the number of variables in any embedded prime of $(P+Q)^2$
is the same as such a number for $P^2$ and $Q^2$.
However,
by \cite[Theorem 1.5]{SW},
$\sqrt{p_i} + \sqrt{q_j}$
are associated to higher powers of $P+Q$,
so that we can get an increase in the number of variables
appearing in an associated prime of higher powers of $P+Q$.

\defin
$(S_z, d_z, S_u, d_u, \varphi)$ is called an {\bf $\bf A$-spreading}
if the following conditions are satisifed:
\ritem
$S_z$ and $S_u$ are $A$-algebras.
\ritem
$d_z, d_u$ are gradings on $S_z, S_u$, respectively,
with degrees in possibly distinct commutative monoids $M_z, M_u$
that are submonoids of free $\bbZ$-modules.
\ritem
$\varphi : S_u \to S_z$ is an $A$-algebra homomorphism
that takes homogeneous components to homogeneous components;
by abuse of notation we write $\varphi : M_u \to M_z$
so that for any homogeneous $f \in S_u$,
$\varphi(d_u(f)) = d_z(\varphi(f))$.
\ritem
If $f$ and $g$ are homogeneous elements in $S_u$,
then the degree of $\varphi(f g)$ is the sum of the degrees
of $\varphi(f)$ and $\varphi(g)$.
In other words,
$d_z (\varphi (f g)) = d_z ( \varphi (f)) + d_z ( \varphi(g))$.
With the abuse of notation from~(3),
this is saying that
$d_z \circ \varphi$ is a monoid homomorphism.

Given such an $A$-spreading,
the {\bf spreading of a homogeneous ideal} $I$ in $S_z$
is
$$
\spread(I) =
\left(f \in S_u: \hbox{$f$ is homogeneous and } \varphi(f) \in I \right).
$$
\endb

\example[exrunning]%
Let $A = k[a,b]$ where $k$ is a field
and $a$ and $b$ are variables over~$k$.
Let $S_z = A[c]$,
$S_u = A[c_1, \ldots, c_n]$,
$d_z$ and $d_u$ trivial on $A$,
$d_z(c) = 1$,
and $\varphi : S_u \to S_z$ given by $\varphi(c_i) = c$ for all $i$.
Let $J = (a^2b^2c, b^4, ab^3, a^3b, a^4) \subseteq S_z$.
\ritem
Let $d_u(c_j) = 1$ for all $j$.
Then 
$\spread(J) = 
(a^2b^2c_1, b^4, ab^3, a^3b, a^4,c_1 - c_2, \ldots, c_1 - c_n)$.
\ritem
Let $d_u(c_j) = e_j$,
where $e_j = (0, \ldots, 0, 1, 0, \ldots, 0)$
has $1$ in the $j$th entry.
Then
$\spread(J) = 
(a^2b^2c_1, \ldots, a^2b^2c_n, b^4, ab^3, a^3b, a^4)$.

Note that $\varphi(c_1 - c_2) = 0 \in J$ regardless of the grading,
but that $c_1 - c_2$ is not in $\spread(J)$ when it is not homogeneous.


\remarks[rmkspread]%
Let $I \subseteq J$ be homogeneous ideals in $S_z$.
\ritem
$\varphi(\spread(I)) \subseteq I$.
Equality need not hold: say if $S_z = A[z]$, $S_u = A[u_1, \ldots, u_n]$
are polynomial rings over $A$
and $\varphi(u_i) = z^2$,
then $\spread(z) = (u_1, \ldots, u_n)$
and
$\varphi(\spread(z)) = \varphi(u_1, \ldots, u_n) = (z)^2 \subsetneq (z)$.
\ritem
If $I \subseteq J$, then $\spread(I) \subseteq \spread(J)$.
\ritem
For all ideals $I, J$ in $S_z$,
$\spread(I) \spread(J) \subseteq \spread(I J)$.

In general,
$\spread(I) \spread(J)$ need not equal $\spread(I J)$
and $\spread(I^r)$ need not equal $\left(\spread(I)\right)^r$.
For example,
let $S_z = A[z]$, $S_u = A[u_1, \ldots, u_n]$
and $\varphi$ take $u_i$ to $z$.
First suppose that $d_u$ is the trivial (zero) grading.
Then $\spread(z^r) = (u_1^r, u_1 - u_2, \ldots, u_1-u_n)$.
The element
$u_1 - u_2$ is thus not in any higher power of $\spread(z)$
but it is in $\spread(z^r)$ for every positive integer $r$.
If instead
we choose the grading $d_u(u_i) = e_i$ for all $i$ and change~$\varphi$
to $\varphi(u_i) = z^2$,
then $\spread(z) = \spread(z^2) = (u_1, \ldots, u_n)$,
so that spreading again does not commute with powers.
Here is another example:
Let $S_z = A[z_1, z_2, \ldots, z_m]$ for some $m \ge 2$,
$S_u = A[u_{ij}: i = 1, \ldots, m; j = 1, \ldots, u_n]$,
$\varphi(u_{ij}) = z_i$
and $d_u(u_{ij}) =~e_j$.
The spreading does not commute with powers
as the spreading of every ideal contains $u_{11} u_{22} - u_{12} u_{21}$.

\lemma[lmspreadcommute]%
Spreading commutes with radicals.
The spreading of a radical (resp.\ prime, primary) ideal is
radical (resp.\ prime, primary).
\endb

\proof
Let $I$ be a homogeneous ideal in $S_z$.
Certainly $\sqrt I$, $\spread(I)$ are homogeneous.

First suppose that $I$ is a radical ideal.
Let $f \in S_u$ be in the radical of $\spread(I)$.
We want to prove that $f \in \spread(I)$.
It suffices to prove the result in case $f$ is homogeneous.
For some large $N$,
$f^N \in \spread(I)$,
so $(\varphi(f))^N = \varphi(f^N) \in I$.
Since $I$ is radical,
it follows that $\varphi(f) \in I$
and so $f \in \spread(I)$.
This proves that the spreading of a radical ideal is radical.

In general,
$\spread(I) \subseteq \spread(\sqrt I)$,
so that the radical of $\spread(I)$
is contained in the radical ideal $\spread(\sqrt I)$.
We next prove the opposite inclusion.
If $f$ in $\spread(\sqrt I)$ is homogeneous,
then $\varphi(f) \in \sqrt I$,
so that for some large integer $N$,
$\varphi(f^N) = \left(\varphi(f)\right)^N \in~I$.
Hence $f^N \in \spread(I)$.
This proves that the radical of $\spread(I)$ equals $\spread(\sqrt I)$,
so that spreading commutes with radicals.

Now let $f, g \in S_u$ be homogeneous such that $f g \in \spread(I)$.
Then $\varphi(f) \varphi(g) = \varphi(f g) \in~I$.

If $I$ is prime,
then either $\varphi(f)$ or $\varphi(g)$ is in $I$,
so that either $f$ or $g$ is in $\spread(I)$.
Thus by homogeneity of ideals, $\spread(I)$ is a prime ideal.

If $I$ is primary and $f \not \in \spread(I)$,
then $\varphi(f) \not \in I$,
so that $\varphi(g) \in \sqrt I$.
But then $g \in \spread(\sqrt I) = \sqrt {\spread(I)}$.
Thus $\spread(I)$ is primary.
\qed

\lemma[lmspreadpd]%
Spreading commutes with intersections of homogeneous ideals.
In particular,
it commutes with (homogeneous) primary decompositions.

Surjective spreading takes irredundant primary decompositions
to irredundant primary decompositions.
\endb

\proof
Let $I = q_1 \cap \cdots \cap q_r$ be an intersection of homogeneous ideals
in $S_z$.
Then $\spread(I) \subseteq \spread(q_1) \cap \cdots \cap \spread(q_r)$.
To prove equality,
let $f$ be a homogeneous element in 
$\spread(q_1) \cap \cdots \cap \spread(q_r)$.
Then $\varphi(f) \in
\varphi\left(\spread(q_1)\right) \cap \cdots \cap
\varphi\left(\spread(q_r)\right)$
$\subseteq q_1 \cap \cdots \cap q_r = I$,
so that $f \in \spread(I)$.
This proves that spreading commutes with intersections.

By \ref{lmspreadcommute},
if $q_i$ is primary,
so is $\spread(q_i)$,
so that spreading commutes with primary decompositions.

Suppose that $\varphi$ is surjective
and that $\spread(I) = \spread(q_2) \cap \cdots \cap \spread(q_r)$.
Let $a \in q_2 \cap \cdots \cap q_r$.
By surjectivity there exists $b \in \spread(q_2 \cap \cdots \cap q_r)$
such that $a = \varphi(b)$.
But then $b \in \spread(q_2) \cap \cdots \cap \spread(q_r) = \spread(I)$,
so that $a = \varphi(b) \in I$,
which says that $q_1$ was an irredundant primary component of $I$.
This proves the lemma.
\qed

The following is a special case of spreading:

\prop[propincreasevarseasy]%
Let $z, z_1, \ldots, z_n$ be variables over a Noetherian ring $A$,
let $I$ be an ideal in $A[z]$
and let $J = I + (z - z_1, \ldots, z - z_n)$
be an ideal in $R = A[z,z_1, \ldots, z_n]$.
Then the following hold for a positive integer $e$.
\ritem
If an associated prime ideal $P$ of~$I^e$ contains~$z$,
then an associated prime ideal of~$J^e$ contains $P$ and
$n+1$ variables $z, z_1, \ldots, z_n$.
\ritem
If $z$ is not a zerodivisor on~$A[z]/I^e$,
then $J^e : z \subseteq I^e + (z - z_1, \ldots, z - z_n)$.
\endb

\proof
Set $S_z = A[z]$,
$S_u = R$,
$d_u$ and $d_z$ both zero on $A$ and with value $1$ on the variables.
Let $\varphi : S_u \to S_z$ be the $A[z]$-algebra homomorphism
that maps $z_i$ to $z$.
Then $(S_z, d_z, S_u, d_u, \varphi)$ is a surjective $A$-spreading,
and $J = \spread(I)$.

Let $P$ be an associated prime ideal of $I^e$ such that $z \in P$.
By a characterization of associated primes,
$P = I^e : a$ for some $a \in S_z \setminus I^e$.
If $a \in J^e$
then $a = \varphi(a) \in \varphi(J^e) = I^e$,
which is a contradiction.
So $J^e : a$ is a proper ideal
and so necessarily contained in some associated prime $Q$ of $J^e$.
This $Q$ contains $z$ and also $z- z_1, \ldots,z- z_n$,
hence it contains $z, z_1, \ldots, z_n$.
(Note that \ref{lmspreadpd} proves~(1) in case $e = 1$
but not in case $e > 1$.)

Suppose that $z$ is not a zerodivisor on~$S_z/I^e$.
Let $b \in S_u$ such that $z b \in J^e$.
Then $z \varphi(b) = \varphi(z b) \in \varphi(J^e) = I^e$.
By assumption then $\varphi(b) \in I^e$,
so that $b \in J^e + (z-z_1,\ldots, z-z_n)$
$= I^e + (z-z_1,\ldots, z-z_n)$.
\qed

\definition
An $A$-spreading $(S_z, d_z, S_u, d_u, \varphi)$ is called {\bf full}
if
for any $d_u$-degree $\underline a$,
$\varphi$~restricted to $\left(S_u\right)_{\underline a}$
is injective into $\left(S_z\right)_{\varphi(\underline a)}$.
\endb

\example[exfull]%
The following is an example of a non-injective full spreading. 
Let $m$, $v_1, \ldots, v_m$ be positive integers,
$S_z = A[z_1, \ldots, z_m]$
and $S_u = A[u_{ij}: i = 1, \ldots, m; j = 1, \ldots, v_i]$.
Let $\varphi(u_{ij}) = z_i$ and $\varphi|_A = id$.
Let $d_z$ be the monomial $\bbN^m$-grading on $S_z$ with $d_z(A) = 0$
and $d_z(z_i) = e_i$,
and let $d_u$ be the monomial $\bbN^N$-grading on $S_u$ with
$d_u(A) = 0$
and $d_u(u_{ij}) = e_{v_1 + v_2 + \cdots + v_{i-1} + j}$.
The verification of the full property is straightforward
since every homogeneous element in $S_z$ (resp.\ $S_u$)
is of the form $a m$
for some $a \in A$ and some monomial $m$ in the $z_i$
(resp.\ in the $u_{ij}$).
%


Recall that for an ideal $J$ in a ring $A$,
its {\bf Rees algebra} is $A[Jt]$,
its {\bf extended Rees algebra} is $A[Jt, t^{-1}]$,
and its {\bf Rees-like algebra} is $A[Jt, t^2]$,
where $t$ is a variable over~$A$.
In the rest of this section we use the
presenting ideal of a Rees algebra of an ideal~$J$
to construct the prime ideal
presenting the same type of Rees algebra of a full spreading of~$J$.


\vskip\parskip
\noindent{\bf Set-up:}
Let $(S_z, d_z, S_u, d_u, \varphi)$ be an $A$-spreading of Noetherian rings.
Let $a_1, \ldots, a_m$ be homogeneous elements in $S_z$
and let $J = (a_1, \ldots, a_m)$.
Let $\spread(J)$ have $m'$ homogeneous generators.
Let $t, T$, $\underline Z = Z_1, \ldots, Z_m$, $U_1, \ldots, U_{m'}$
be variables over~$S_z$ and~$S_u$.

Let $\widetilde \varphi:
S_u[t, t^{-1}] \to S_z[t, t^{-1}]$
be the $A$-algebra homomorphism that agrees with~$\varphi$ on~$S_u$
and such that $\varphi(t^{\pm 1}) = t^{\pm 1}$.
The grading on $S_z[t,t^{-1}]$ is as follows:
the degree of a homogeneous element $s \in S_z \subseteq S_z[t, t^{-1}]$
is $(d_z(s), 0)$
and the degree of~$t^{\pm 1}$ is~$(0,\pm 1)$.
The grading on $S_u[t,t^{-1}]$ is defined similarly.

All three types of Rees algebras of $J$ are subrings
of $S_z[t, t^{-1}]$.
Let
$\psi_z : S_z[\underline Z, T] \to S_z[t,t^{-1}]$
be the $S_z$-algebra homomorphism
which takes $Z_i$ to $a_i t$,
and which takes $T$ to one of $0$, $t^{-1}$, $t^2$,
depending on whether we are using
Rees algebra, extended Rees algebra, or the Rees-like algebra.
The image of $\psi_z$ is the chosen type of the Rees algebra of~$J$.
To make this map graded,
we impose the grading $d_Z$ on $S_z[\underline Z, T]$ as follows:
$d_Z(s) = (d_z(s), 0)$ for all $s \in S_z$,
$d_Z(Z_i) = (d_z(a_i),1)$ for all $i = 1, \ldots, m$,
and
the degree of $T$ is
$\infty$, $(0,-1)$, or $(0, 2)$ depending on the type of the Rees algebra.
(Instead of $d_Z(T) = \infty$ in case of a Rees algebra
we can simply not adjoin $T$.)

Analogously we define the grading $d_U$ on $S_u[\underline U, T]$,
and we let $\psi_u: S_u[\underline U, T] \to S_u[t, t^{-1}]$
be a graded $S_u$-algebra homomorphism
whose image is a Rees algebra of $\spread(J)$,
this Rees algebra being of the same type as the constructed algebra for~$J$.
With this we have the following commutative diagram:
$$
\matrix{
J \subseteq S_z & \subseteq & S_z[\underline Z, T]
& {\psi_z \atop \longrightarrow}
& \hbox{Rees algebra of $J$} & \subseteq & S_z[t,t^{-1}] \cr
\cr
\ \ \ \ \ \ \ \ \uparrow \varphi 
& & & & & & \ \ \ \ \uparrow \widetilde \varphi \cr
\cr
\spread(J) \subseteq S_u & \subseteq & S_u[\underline U, T]
& {\psi_u\atop \longrightarrow}
& \hbox{Rees algebra of $\spread(J)$} & \subseteq & S_u[t,t^{-1}] \cr
}
$$
Note that $\widetilde \varphi$ takes the Rees algebra of $\spread(J)$
to the Rees algebra of~$J$.
Let $\widehat \varphi$ be the restriction map of $\widetilde \varphi$.
We next define $\Phi: S_u[\underline U, T] \to S_z[\underline Z, T]$
as
$\Phi|_{S_u} = \varphi$,
$\Phi(T) = T$,
and $\Phi(U_i) = f_i$
where $f_i$ is a homogeneous element of $S_z[\underline Z, T]$
such that $\psi_z(f_i) = \widehat \varphi(\psi_u(U_i))$.
Then the following diagram commutes:
$$
\matrix{
J \subseteq S_z & \subseteq & S_z[\underline Z, T]
& {\psi_z \atop \longrightarrow}
& \hbox{Rees algebra of $J$} & \subseteq & S_z[t,t^{-1}] \cr
\cr
\ \ \ \ \ \ \ \ \uparrow \varphi 
& & \uparrow \Phi & & \uparrow \widehat \varphi & & \ \ \ \ \uparrow \widetilde \varphi \cr
\cr
\spread(J) \subseteq S_u & \subseteq & S_u[\underline U, T]
& {\psi_u\atop \longrightarrow}
& \hbox{Rees algebra of $\spread(J)$} & \subseteq & S_u[t,t^{-1}] \cr
}
$$

\thm[thmspreadReescommute]%
$(S_z[\underline Z, T], d_Z, S_u[\underline U, T], d_U, \Phi)$
is an $A$-spreading.
\ritem
If $\varphi$ is surjective,
then there exists a surjective $\Phi$.
\ritem
$\ker(\psi_u) \subseteq \spread(\ker(\psi_z))$.
\ritem
If $\varphi$ is full,
then $\spread(\ker(\psi_z))$ equals $\ker(\psi_u)$.
In other words,
when $\varphi$ is full then
the spreading of the presenting ideal of a type of Rees algebra of $J$
equals the presenting ideal of the same type of Rees algebra
of the spreading of $J$.
\endb

\proof
The set-up makes
$(S_z[\underline Z, T], d_Z, S_u[\underline U, T], d_U, \Phi)$ a spreading.

Suppose that $\varphi$ is surjective.
Then $\widetilde \varphi$ and $\widehat \varphi$ are surjective as well.
We may order a generating set of $\spread(J)$
so that the image under $\varphi$ of the $i$th generator is $a_i$
for $i \le m$.
We then set $\Phi(U_i) = Z_i$ for $i \le m$,
which makes $\Phi$ surjective.


Let $f \in \ker(\psi_u)$.
Since $\psi_u$ is a graded homomorphism,
every homogeneous component of $f$ is in $\ker(\psi_u)$,
so that to prove that $f \in \spread(\ker(\psi_z))$
without loss of generality we may assume that $f$ is homogeneous.
Then
$\psi_z \circ \Phi (f) = \widehat \varphi \circ \psi_u(f) = 0$,
so that $\Phi(f) \in \ker (\psi_z)$.
Thus $f \in \spread(\ker(\psi_z))$.

Now let $\varphi$ be full
and let $f \in \spread(\ker(\psi_z))$ be homogeneous.
Then $\Phi(f) \in \ker(\psi_z)$,
so that $0 = \psi_z \circ \Phi (f) = \widehat \varphi \circ \psi_u(f)$.
Since $\varphi$ is injective on homogeneous components,
so are $\widetilde \varphi$ and $\widehat \varphi$,
so that $\psi_u(f) = 0$.
Thus $f \in \ker(\psi_u)$.
\qed


\section{Examples}[sectexamples]

This section provides a few examples of prime ideals
whose powers have embedded primes that contain variables.
These examples are used
in \ref{sectsplexamples}
as a base for generating prime ideals
whose powers have many associated primes containing many variables.

The first example below treats all powers of a prime ideal,
whereas the second example is about the second power only.

\example[exsimpleP]%
Let $P$ be the kernel of the $k$-algebra homomorphism
$k[x, y, z] \to k[t^3,t^4,t^5]$ taking
$x$ to $t^3$,
$y$ to $t^4$,
$z$ to $t^5$.
Then $P$ is a prime ideal of height two,
it contains no variables,
and by~\cite{Hunpc3},
$(x,y,z)$ is associated to $P^e$ for all $e \ge 2$,
and $\cup_{e=1}^\infty \Ass(k[x,y,z]/P^e) = \{P, (x_1,y,z)\}$.

\prop[propReesP]%
In any polynomial ring in nine variables over an arbitrary field
there exists a binomial prime ideal~$P$ of height~$5$ containing no variables
such that $P^2$ has exactly two embedded associated prime ideals:
one of the two is a monomial ideal generated by eight variables
and the other is the maximal ideal generated by the nine variables.
\endb

\proof
Let $k$ be a field,
$a, b, c, Z_1, \ldots, Z_5, t, T$ variables over $k$,
let
$A = k[a,b,c]$,
$J = (a^2b^2c, b^4, ab^3, a^3b, a^4)A$,
and $R = k[a,b,c,Z_1,Z_2,Z_3,Z_4,Z_5,T]$.
Let $B = A[Jt,t^{-1}] \subseteq A[t,t^{-1}]$
be the extended Rees algebra of $J$.
There is a natural surjective $A$-algebra map $R \to B$,
with $Z_i$ mapping onto $t$ times the $i$th listed generator of $J$,
and with $T$ mapping to $t^{-1}$.
Let $P$ be the kernel of this map.

Then $P$ is a prime ideal that contains no variables.
As the dimension of the Rees algebra is~$4$,
the height of $P$ is~$5$.
Rees algebras of monomial ideals are generated by binomials.

It is easy to verify that the following elements are in~$P$:
$f_1 = a^4 - Z_5 T$,
$f_2 = aZ_2 - bZ_3$,
$f_3 = aZ_4 - bZ_5$,
$f_4 = ab^3 - Z_3 T$,
$f_5 = a^3Z_3 - b^3Z_5$,
$f_{6} = a^4 Z_2 - b^4 Z_5$,
$f_{7} = a^2 Z_1 - b^2 c Z_5$,
$f_{8} = a^2 b^2 c -Z_1 T$,
$f_{9} = Z_1 Z_5 - c Z_4^2$,
$f_{10} = Z_1 Z_2 - c Z_3^2$,
$f_{11} = aZ_1-bc Z_4$,
$f_{12} = Z_1^2 - c^2 Z_3 Z_4$,
$f_{13} = Z_2 T - b^4$,
$f_{14} = Z_2^2 Z_4 - Z_3^3$,
$f_{15} = Z_2^3 Z_5 - Z_3^4$,
$f_{16} = acZ_3-b Z_1$,
$f_{17} = a^3 b -Z_4 T$,
$f_{18} = a^2Z_3 - b^2Z_4$.
We do not claim that these elements generate~$P$.

Consider
$\alpha = a^4 Z_2-a^3 b Z_3-a b^3 Z_4+b^4 Z_5-Z_2 Z_5 T+Z_3 Z_4 T$.
If $\alpha \in P^2$,
then under the lexicographic order $a > b > c > Z_1 > Z_2 > \cdots > Z_5 > T$,
the leading monomial $a^4 Z_2$ must be a product of two leading terms of
elements of $P$.
By the structure of the kernels of monomial maps,
neither $a^3$ nor $Z_2$ can be multiples of leading terms of elements of $P$,
which means that 
$\alpha$ is not in $P^2$.
One can verify that
$$
\eqalignno{
a\alpha &= 
a^5 Z_2-a^4 b Z_3-a^2 b^3 Z_4+a b^4 Z_5-a Z_2 Z_5 T+a Z_3 Z_4 T
= f_1 f_2 - f_3 f_4 \in P^2.
\cr
}
$$
Since $a \not \in P$,
this proves that $P^2$ has an embedded prime ideal
which contains $P^2 : \alpha$ and thus $a$.
If we invert $c$ or set $c = 1$,
then the resulting $\alpha$ is still not in $P^2$
but $a \alpha \in P^2$.
This proves that $c$ is not in the radical of $P^2 : \alpha$,
so that at least one of the embedded prime ideals of $P^2$
does not contain~$c$.

After inverting $a$,
the ideal $(f_1, f_3, f_5, f_{6},f_{7})$
equals
$(a^{-4}Z_5 T - 1,
Z_4 - b Z_5/a,
Z_3 - b^3 Z_5/a^3,
Z_2 - b^4 Z_5/a^4,
Z_1 - b^2c Z_5/a^2)$,
which is a prime ideal in $R_a$ of height~$5$ and contained in $P_a$.
By the height consideration
$P_a$ must equal this five-generated ideal.
Similarly,
after inverting $Z_1$,
$P$ is the complete intersection prime ideal
$(f_{8}, f_{9}, f_{10},f_{11}, f_{12})_{Z_1}
= (T - a^2 b^2 c Z_1^{-1},
Z_5 - c Z_4^2 Z_1^{-1},
Z_2 - c Z_3^2 Z_1^{-1},
a - bc Z_4 Z_1^{-1},
1 - c^2 Z_3 Z_4 Z_1^{-1})$,
after inverting~$Z_2$,
$P$ is the complete intersection prime ideal
$(f_2, f_{10}, f_{13}, f_{14}, f_{15})_{Z_2}
= (a - b Z_3 Z_2^{-1},
Z_1 - c Z_3^2 Z_2^{-1},
T - b^4 Z_2^{-1},
Z_4 - Z_3^3 Z_2^{-2},
Z_5  - Z_3^4 Z_2^{-3})$,
and after inverting~$Z_3$,
$P$ is the complete intersection prime ideal
$(f_2 , f_4, f_{10}, f_{14}, f_{15})_{Z_3}
= (b - a Z_2 Z_3^{-1},
T - a b^3 Z_3^{-1},
c - Z_1 Z_2 Z_3^{-2},
Z_5 - Z_3^4 Z_3^{-3},
1 - Z_2^2 Z_4 Z_3^{-3})$.
Furthermore,
after localization at $T$,
$P_T$ is generated by the variables
$Z_1 - a^2b^2cT^{-1}$,
$Z_2 - b^4T^{-1}$,
$Z_3 - ab^3 T^{-1}$,
$Z_4 - a^3b T^{-1}$,
$Z_5 - a^4 T^{-1}$.
Whenever an ideal in a Cohen-Macaulay ring is generated by a regular
sequence,
its powers have no embedded primes.
So we just proved that $a$, $Z_1$, $Z_2$ and $Z_3, T$
must be contained in every embedded prime of every power of $P$.
By symmetry,
$a,b,Z_1,Z_2,Z_3,Z_4,Z_5, T$
must be contained in all the embedded primes of powers of $P$.

Thus we have proved that $P^2$ has an embedded associated prime ideal
and that each embedded prime ideal contains 
$(a, b, Z_1, Z_2, Z_3, Z_4, Z_5, T)$.
By multihomogeneity of $J$ and of the extended Rees algebra,
then the only possible associated primes of $P^2$ are
$Q_1 = (a, b, Z_1, Z_2, Z_3, Z_4, Z_5, T)$
and $Q_2 = (a, b, c, Z_1, Z_2, Z_3, Z_4, Z_5, T)$.
Since we proved that at least one embedded prime ideal does not contain~$c$,
we get that $Q_1$ is associated to $P^2$.

We next prove that $Q_2$ is also associated to $P^2$.
Consider
$\beta = a^5 Z_3-2 a^3 b^2 Z_4+a^2 b^3 Z_5-a Z_3 Z_5 T+b Z_4^2 T$.
If $\beta \in P^2$,
then
under the lexicographic order $a > b > c > Z_1 > Z_2 > \cdots > Z_5 > T$,
the leading monomial $a^5 Z_3$ must be a product of the leading terms of
two elements of $P$.
By the structure of the kernels of monomial maps,
neither $a^3$ nor $a Z_3$ can be multiples of leading terms of elements of $P$,
which means that $\beta \not \in P^2$.
However,
$$
c \beta = 
a^5 c Z_3-2 a^3 b^2 c Z_4+a^2 b^3 c Z_5-a c Z_3 Z_5 T+b c Z_4^2 T
= f_1 f_{16} 
+ f_{11} f_{17} 
- f_{8} f_{18}.
$$
This proves that $Q_2$ is associated to $P^2$.
\qed

The proof shows more:

\prop[propReesPmax]%
Let $k$ be an arbitrary field,
let $a, b, c, Z_1, \ldots, Z_5, T$ be variables over $k$,
let $R = k[a, b, c, Z_1, \ldots, Z_5, T]$,
and let $P$ be the kernel of the $k[a,b,c]$-algebra surjection 
from~$R$ to the extended Rees algebra of
the monomial ideal $(a^2b^2c, b^4, ab^3, a^3b, a^4) \subseteq k[a,b,c]$,
where $Z_i$ maps to the $i$th listed generator of $J$
and $T$ maps to $t^{-1}$.
Then there exist $\alpha, \beta \in k[a,b, Z_2,Z_3, Z_4, Z_5, T]$
such that $(a,b,Z_1, \ldots, Z_5, T)$
is the radical of $P^2 : \alpha$
and such that the maximal ideal $(a,b,c,Z_1, \ldots, Z_5, T)$
is the radical of $P^2 : \beta$.
(We emphasize:
no variables $c$ and $Z_1$ appear in $\alpha$ and $\beta$.)
\endb

\proof
\ref{propReesP} shows that $P^2 = P^{(2)} \cap J_1 \cap J_2$,
where $J_1$ is primary to $(a, b, Z_1, \ldots, Z_5, T)$
and $J_2$ is primary to $(a, b, c, Z_1, \ldots, Z_5, T)$.

We take $\alpha, \beta$ as in the proof of \ref{propReesP}.
Since $a \alpha \in P^2$
and since $a$ is a non-zerodivisor on $R/P^{(2)}$,
it follows that $\alpha \in P^{(2)}$.
Thus $P^2 : \alpha = (J_1 : \alpha) \cap (J_2 : \alpha) \not = R$.
Since $c$ is not in the radical of $P^2 : \alpha$
(see the proof of \ref{propReesP}),
necessarily $J_1 : \alpha \not = R$.
Thus the radical of $P^2 : \alpha$ is $(a,b,Z_1, \ldots, Z_5, T)$.

Since $c \beta \in P^2$
and since $c$ is a non-zerodivisor on $R/P^{(2)}$ and on $R/J_1$,
it follows that $\beta \in P^{(2)} \cap J_1$.
Thus $P^2 : \beta = J_2 : \beta \not = R$
must be primary to the maximal ideal.
\qed

\remark
Work similar to that in the proof of \ref{propReesP} shows that
if $P$ is the presenting ideal of the Rees algebra $A[Jt]$ of~$J$,
then $P$ is a prime ideal of height~$4$ in a polynomial ring in eight variables,
that
$(a, b, Z_1, Z_2, Z_3, Z_4, Z_5)$ is associated to $P^2$
and that the only other candidate
for an embedded associated prime of $P^2$
is $(a, b, c, Z_1, Z_2, Z_3, Z_4, Z_5)$.
Macaulay2~\cite{GS} computes that the latter ideal is not associated to $P^2$.

If $P$ is the presenting ideal of the Rees-like algebra $A[Jt,t^2]$ of~$J$,
then we can similarly show
that $P$ is a prime ideal of height~$5$ in a polynomial ring in nine variables,
that either 
$(a, b, Z_1, Z_2, Z_3, Z_4, Z_5)$
or $(a, b, Z_1, Z_2, Z_3, Z_4, Z_5, T)$
is associated to $P^2$
and that the only other candidates
for embedded associated primes of $P^2$
are $(a, b, c, Z_1, Z_2, Z_3, Z_4, Z_5)$
and $(a, b, c, Z_1, Z_2, Z_3, Z_4, Z_5, T)$.
Macaulay2~\cite{GS} computes that only the first and the third prime ideals
on this list are associated to $P^2$.

\section{Prime ideals whose powers have many associated prime ideals}[sectsplexamples]%

We exploit splitting and spreading to generate prime ideals whose
specific powers
have arbitrarily many associated prime ideals.

\thm[thmmainone]
Let $k$ be an arbitrary field,
and let $m \ge 3$ and $v_1, \ldots, v_m$ be any positive integers.
Then there exists a polynomial ring $R$ in $\sum v_i$ variables over $k$
with an $m$-generated prime ideal~$P$ of height $m-1$
with generators of degree at most~$3$
(or with quasi-homogeneous generators of degree at most~$10$)
such that for all integers $e \ge 2$,
$P^e$ has $\prod v_i$ embedded primes,
all of which have height $m$.

In case all $v_i$ equal $v$,
this says that there exists a polynomial ring in $m v$ variables
over a field $k$
with an $m$-generated prime ideal~$P$ of height $m-1$
such that $P^e$ has $v^m$ embedded primes if $e \ge 2$.
The number $v^m = \left({\root v \of v}\right)^{mv}$
is exponential in the number $mv$ of variables
if we think of $v$ as fixed.
\endb

\proof
Let $x$, $y$, $z, z_1, \ldots, z_{m-3}$ be variables over $k$.
Let $I = (x^3-yz, y^2-xz, z^2-x^2y) \in k[x,y,z]$
be the prime ideal as in \ref{exsimpleP}.
Set $R_m = k[x,y,z, z_1, \ldots z_{m-3}]$
and $I_m = (x^3-yz, y^2-xz, z^2-x^2y,z_1-z, \ldots, z_{m-3} - z)$.
Then $R_m$ is a polynomial ring in $m$ variables over $k$
and $I_m$ is a prime ideal in $R_m$ of height $m-1$.
By \ref{propincreasevarseasy},
$I_m$ is a spreading of~$I$.
By \ref{exsimpleP},
for all integers~$e \ge 2$,
$I^e$ has exactly one embedded prime ideal, namely $(x,y,z)$.
Then by \ref{propincreasevarseasy},
the maximal ideal $J = (x,y,z,z_1, \ldots, z_{m-3})$ in $R_m$
is an associated prime ideal of $I_m^e$.
Suppose that $Q$ is associated to $I_m^e$.
Since $I_m$ has height $m-1$,
the height of $Q$ is either $m-1$ or $m$,
so that by quasi-homogeneity of $I_m$
we get that either $Q = I_m$ or $Q = J$.
This proves that for all integers $e \ge 2$,
the set of associated primes of $I_m^e$
consists of $I_m$ and $J$.

Let $\varphi : R_m \to A$ be the splitting
which for $i = 1, \ldots, m$
splits the $i$th variable into $v_i$ variables.
Set $P = \varphi(I_m)A$.
Then by \ref{thmpdgensplitmorevars},
$A$ is a polynomial ring in $\sum_i v_i$ variables over $k$,
$P$ is a prime ideal of height $m-1$,
and $P^e$ has $1 + v_1 \cdots v_m$ primary components for all $e \ge 2$.
Each of the embedded components is the splitting of $J$
and is thus generated by $m$ variables.

The number of generators of $P$
is at most the number of generators $m$ of $I_m$.
Since the height of $P$ is $m-1$,
the number of generators is at least $m-1$,
and since higher powers of $P$ have embedded primes,
$P$ cannot be a complete intersection,
so that the number of generators of $P$ is exactly $m$.
\qed

\cor
Let $M$ and $d$ be positive integers
and let $E(M,d)$ the Ananyan-Hochster bound from~\cite{AH},
namely the constant such that for any ideal~$I$
in any polynomial ring over a field,
if $I$ has at most $M$ generators all of which have degree at most~$d$,
then $I$ has at most $E(M,d)$ associated primes.
Then $E(M,d) \ge 1 + 3^{\sqrt{2M}-1}$.
\endb

\proof
Let $m \ge 3$, $v = 3$, $P$ be as in \ref{thmmainone}.
Then the number $M$ of generators of~$P^2$ is at most
${m(m+1) \over 2} \le {(m+1)^2 \over 2}$,
and the number of associated primes of $P^2$ is
$$
1 + 3^m \ge 1 + 3^{\sqrt{2M} - 1}.
\eqed
$$

\thm[thmmainmain]%
For every field $k$,
for every odd integer $m = 2n + 7 \ge 9$
and for all positive integers $v_1, \ldots, v_m$
there exists a polynomial ring in $\sum_{i = 1}^m v_i$ variables over~$k$
with a prime ideal $P$ of height ${m+1 \over 2} = n + 4$
such that $P^2$ has at least
$\prod_{i = 1}^m v_i + v_1 v_2 \prod_{i = n+3}^m v_i$ embedded primes.

If all $v_i$ equal $v$,
this says that there exists a prime ideal $P$
of height ${m+1 \over 2} = n + 4$
in a polynomial ring in $mv$ variables
such that $P^2$ has at least $v^m + v^{{m+7 \over 2}}$
$= v^{2n+7} + v^{n+7}$ embedded primes.
The number
$v^m = \left({\root v \of v}\right)^{mv}$
is exponential in the number $mv$ of variables
if we think of $v$ as fixed.
\endb

\proof
We use the set-up as in \ref{exrunning}~(2):
$A = k[a,b]$,
$S_z = A[c]$,
$S_u = A[c_1, \ldots, c_n]$,
$J = (a^2b^2c, b^4, ab^3, a^3b, a^4) \subseteq S_z$,
$d_z$ and $d_u$ trivial gradings on $A$,
and $d_z(c) = 1$,
$d_u(c_j) =~e_j$,
where $e_j = (0, \ldots, 0, 1, 0, \ldots, 0)$ has $1$ in the $j$th entry.
Then the $A$-module homomorphism $\varphi : S_u \to S_z$
taking $c_j$ to $c$ is a spreading,
and it is surjective and full.
Furthermore,
$\spread(J) = (a^2b^2c_1, \ldots, a^2b^2c_n, b^4, ab^3, a^3b, a^4)$.

The extended Rees algebra of $J$
is a natural homomorphic image of the polynomial ring
$S_Z = S_z[Z_1, \ldots, Z_5,T]$,
and the extended Rees algebra of $\spread(J)$
is a natural homomorphic image of the polynomial ring
$S_U = S_u[U_1, \ldots, U_n, U'_2, U'_3, U'_4,U'_5,T]$.
Let $P_Z$ (resp.\ $P_U$) be the kernels of these homomorphisms.
By \ref{thmspreadReescommute},
there exists a surjective spreading
$(S_Z,d_Z,S_U,d_U,\Phi)$
such that $\spread(P_Z) = P_U$.
We may take $\varphi(U_i) = Z_1$
and $\varphi(U'_i) = Z_i$ for all $i$.
Without loss of generality we may identify each $U'_i$ with $Z_i$,
so that in the sequel we write
$S_U = S_u[U_1, \ldots, U_n, Z_2, Z_3, Z_4,Z_5,T]$.

Note that $S_U$ is a polynomial ring in $m = 2n+7$ variables over $k$.
Since $S_u$ has dimension $n+2$,
the extended Rees algebra of $\spread(J)$ has dimension $n+3$,
so that the height of $P_U$ is $(2n+7) - (n+3) = n+4 = {m+1 \over 2}$.

Let $M_Z$ (resp.\ $M_U$)
be the maximal homogeneous ideal in $S_Z$ (resp.\ $S_U$).
We show that $M_U$ is associated to $P_U^2$.
Note that $P_Z$ is as in \ref{propReesPmax},
and so $M_Z$ is associated to $P_Z^2$
and $M_Z$ is the radical of $P_Z^2 : f$
for some $f \in k[a,b,Z_2, Z_3, Z_4, Z_5,T]$.
Since $\Phi(f) = f \not \in P_Z^2 = \Phi(P_U^2)$,
it follows that $f \not \in P_U^2$.
Let $x$ be any variable in~$S_U$.
Then for some large integer $p$,
$\Phi(x^p f) = \Phi(x)^p f \in P_Z^2 = \Phi(P_U^2)$.
We use a lexicographic order that places all $c_1, \ldots, c_n$ and
$U_1, \ldots, U_n$ at the top of the order,
and if $x$ is one of these variables,
then $x$ is the least in the order.
Then by Gr\"obner bases theory
there exists an equation
that writes $\Phi(x^p f)$ as an element of $P_Z^2$
using only the variables $\Phi(x),a,b,Z_2, Z_3, Z_4, Z_5,T$.
If in that equation we replace each occurrence of $\Phi(x)$ with $x$,
then by the definition of spreading we get that $x^p f \in P_U^2$.
This proves that $M_U$ is associated to $P_U^2$.

Set $Q_Z = (a,b,Z_1, \ldots, Z_5, T)$,
and $Q'_U = (a,b,U_1, \ldots, U_n, Z_2, \ldots, Z_5, T)$.
By \ref{propReesPmax},
$Q_Z$ is associated to $P_Z^2$
and $Q_Z$ is the radical of $P_Z^2 : g$
for some $g \in k[a,b,Z_2, Z_3, Z_4, Z_5,T]$.
As in the proof in the previous paragraph,
$g \not \in P_U^2$
and a power of $Q'_U$ is contained in $P_U^2 : g$.
If a power of $c_i$ is in $P_U^2 : g$
then a power of $c = \Phi(c_i)$ is in $P_Z^2 : g$,
which is a contradiction.
Thus no power of $c_i$ is in $P_U^2 : g$.
Thus there exists a prime ideal $Q_U$ associated to $P_U^2$
that contains $Q'_U$ and that is different from~$M_U$.
Thus $Q_U$ contains at least the $n+7 = {m-7 \over 2}$ variables
from $Q'_U$.

If we split the $i$th variable in $S_U$
into a product of
$v_i$ distinct new variables,
as $i$ varies from $1$ to $m$,
then by \ref{thmpdgensplitmorevars},
the image $P$ of $P_U$ is a prime ideal of height ${m+1 \over 2}$,
and $P^2$ has at least $\prod_{i = 1}^m v_i$ embedded primes
coming from $M_U$
and at least
$v_1 v_2 \prod_{i = n+3}^m v_i$ embedded primes
coming from $Q_U$.
\qed

\remark
The proof shows that one of the associated prime ideals
of $P_U^2$ is the maximal ideal of $S_U$
and that another associated prime ideal has at least $n + 7$ variables.
We do not determine the exact number of variables
in this second associated prime.

\examples
All rings below are polynomial rings over an arbitrary field.
\ritem
By \ref{thmmainone}
there exists a prime ideal~$P$ of height~$4$
in $2 + 2 + 5 + 5 + 5 = 19$ variables
such that $P^e$ has exactly $500 = 2^2 \cdot 5^3$ embedded primes
for all $e \ge 2$.
\ritem
By \ref{thmmainmain}
there exists a prime ideal~$P$ of height~$5$
in $2 + 2 + 1 + 1 + 1 + 1 + 5 + 5 + 5 = 23$ variables
such that $P^2$ has at least
$2^2 \cdot 1^4 \cdot 5^3 + 2^2 \cdot 1^3 \cdot 5^3 = 1000$ embedded primes.
\ritem
By \ref{thmmainone},
there exists a prime ideal $P$ of height~$5$
in $6 \cdot 3 = 18$ variables
such that $P^e$ has exactly $3^6 = 729$ embedded primes for all $e \ge 2$.
\ritem
By \ref{thmmainmain},
there exists a prime ideal $P$ of height~$5$
in $9 \cdot 2 = 18$ variables
such that $P^2$ has at least $2^9 + 2^8 = 768$ embedded primes.
\ritem
By \ref{thmmainone},
there exists a prime ideal $P$ of height $10$
in $11 \cdot 2 = 22$ variables
such that $P^e$ has exactly $2^{11} = 2048$ embedded primes for all $e \ge 2$.
\ritem
By \ref{thmmainmain},
there exists a prime ideal $P$ of height~$6$
in $11 \cdot 2 = 22$ variables
such that $P^2$ has at least $2^{11} + 2^9 = 2560$ embedded primes.

\bigskip\bigskip
\leftline{\bf References}
\bigskip

\font\eightrm=cmr8 \def\rm{\fam0\eightrm}
\font\fiverm=cmr5 \def\rm{\fam0\eightrm}
\font\eightit=cmti8 \def\it{\fam\itfam\eightit}
\font\eightbf=cmbx8 \def\bf{\fam\bffam\eightbf}
\font\eighttt=cmtt8 \def\tt{\fam\ttfam\eighttt}
\textfont0=\eightrm \scriptfont0=\fiverm
\rm
\baselineskip=9.9pt
\parindent=3.6em
\setbox1=\hbox{[999]} 
\newdimen\labelwidth \labelwidth=\wd1 \advance\labelwidth by 2.5em
\ifnumberbibl\advance\labelwidth by -2em\fi

\thmno=0

\bitem{AH}
T. Ananyan and M. Hochster,
Small subalgebras of polynomial rings and Stillman's conjecture,
{\tt arXiv:math.AC/1610.09268}.

\bitem{BS}
D.\ Bayer and M.\ Stillman,
On the complexity of computing syzygies,
{\it J.\ Symbolic Comput.}, {\bf 6} (1988), 135-147.

\bitem{Eisenbud}
D. Eisenbud,
{\bkt Commutative Algebra with a View toward Algebraic Geometry},
Springer-Verlag, 1994.

\bitem{GS}
D. Grayson and M. Stillman,
Macaulay2, a software system for research in algebraic geometry,
available at {\tt http://www.math.uiuc.edu/Macaulay2}.

\bitem{H}
G. Hermann,
Die Frage der endlich vielen Schritte in der Theorie der Polynomideale,
{\it Math. Ann.} {\bf 95} (1926), 736--788.

\bitem{Hunpc3}
C. Huneke,
The primary components of and integral closures of ideals in 3-dimensional
regular local rings,
{\it Math. Ann.} {\bf 275} (1986), 617--635.

\bitem{K}
G. A. Kirkup,
Minimal primes over permanental ideals,
{\it Trans. Amer. Math. Soc.} {\bf 360} (2008), 3751--3770.

\bitem{LS}
R. C. Laubenbacher and I. Swanson,
Permanental ideals,
{\it J. Symbolic Comput.} {\bf 30} (2000), 195--205.

\bitem{Mats}
H. Matsumura,
{\bkt Commutative Ring Theory},
Cambridge University Press, 1986.

\bitem{MM}
E.\ Mayr and A.\ Meyer,
The complexity of the word problems for commutative semigroups
and polynomial ideals,
{\it Adv.\ Math.}, {\bf 46} (1982), 305-329.

\bitem{McP}
J. McCullough and I. Peeva,
Counterexamples to the regularity conjecture,
{\it J. Amer. Math. Soc.} {\bf 31} (2018), 473--496. 

\bitem{PS}
J. Porcino and I. Swanson,
2 x 2 permanental ideals of hypermatrices,
{\it Comm. Alg.} {\bf 43} (2015), 84--101.

\bitem{Sei74}
A. Seidenberg,
Constructions in algebra,
{\it Trans. Amer. Math. Soc.} {\bf 197} (1974), 273--313.

\bitem{S2}
I. Swanson,
The minimal components of the Mayr-Meyer ideals,
{\it J.\ Algebra} {\bf 267} (2003), 127--155.

\bitem{S3}
I. Swanson,
On the embedded primes of the Mayr-Meyer ideals,
{\it J. Algebra} {\bf 275} (2004), 143-190.

\bitem{SW}
I. Swanson and R. M. Walker,
Tensor-Multinomial Sums of Ideals:
Primary Decompositions and Persistence of Associated Primes,
{\tt arXiv:1806.03545}.

\bitem{V}
L. G. Valiant,
The complexity of computing the permanent,
{\it Theoret. Comp. Sci.} {\bf 8} (1979), 189--201.

\bitem{vdDS}
L. van den Dries and K. Schmidt,
Bounds in the theory of polynomial rings over fields.
A nonstandard approach,
{\it Invent. Math.} {\bf 76} (1984), 77--91.

\bitem{W}
R. M. Walker,
Uniform symbolic topologies via multinomial expansions,
{\tt arXiv:math.AC/1703.04530}.

\end